\documentclass[12pt,reqno]{amsart}

\usepackage{latexsym}

\newcommand{\leqsim}{\begin{array}{c}<\\[-8pt] \sim\end{array}}
\newcommand{\boxb}{\raisebox{-.8pt}{$\Box$}_b}
\newcommand{\sm}{\setminus}
\newcommand{\szego}{Szeg\"o }
\newcommand{\nhat}{\raisebox{2pt}{$\wh{\ }$}}
\newcommand{\Si}{\Sigma}
\newcommand{\zetaone}{\zeta^{(1)}}
\newcommand{\zetatwo}{\zeta^{(2)}}
\newcommand{\inv}{^{-1}}
\newcommand{\kahler}{K\"ahler }
\newcommand{\sqrtn}{\sqrt{N}}
\newcommand{\wt}{\widetilde}
\newcommand{\wh}{\widehat}
\newcommand{\PP}{{\mathbb P}}
\newcommand{\R}{{\mathbb R}}
\newcommand{\C}{{\mathbb C}}

\newcommand{\CP}{\C\PP}
\renewcommand{\d}{\partial}
\newcommand{\dbar}{\bar\partial}
\newcommand{\ddbar}{\partial\dbar}

\renewcommand{\H}{{\mathbf H}}
\newcommand{\zb}{{\mathbf Z}}
\newcommand{\half}{{\frac{1}{2}}}

\newcommand{\diag}{{\operatorname{diag}}}

\renewcommand{\phi}{\varphi}
\newcommand{\eqd}{\buildrel {\operatorname{def}}\over =}

\newcommand{\ccal}{\mathcal{C}}

\newcommand{\gcal}{\mathcal{G}}
\newcommand{\hcal}{\mathcal{H}}
\newcommand{\ical}{\mathcal{I}}
\newcommand{\lcal}{\mathcal{L}}

\newcommand{\jcal}{\mathcal{J}}
\newcommand{\al}{\alpha}
\newcommand{\be}{\beta}
\newcommand{\ga}{\gamma}

\newcommand{\La}{\Lambda}
\newcommand{\la}{\lambda}
\newcommand{\ep}{\varepsilon}
\newcommand{\de}{\delta}

\newcommand{\om}{\omega}
\newcommand{\Om}{\Omega}
\newcommand{\di}{\displaystyle}
\newtheorem{theo}{{\sc Theorem}}[section]
\newtheorem{maintheo}{{\sc Theorem}}

\newtheorem{lem}[theo]{{\sc Lemma}}
\newtheorem{prop}[theo]{{\sc Proposition}}

\newenvironment{rem}{\medskip\noindent{\it Remark:\/} }{\medskip}
\newenvironment{defin}{\medskip\noindent{\it Definition:\/} }{\medskip}

\title[Asymptotics of almost holomorphic sections on symplectic manifolds]
{Asymptotics of almost holomorphic sections of ample line bundles on symplectic
manifolds}

\author{Bernard Shiffman}
\author{Steve Zelditch}
\address{Department of Mathematics, Johns Hopkins University, Baltimore,
MD
21218, USA}
\email{shiffman@math.jhu.edu, zelditch@math.jhu.edu}

\thanks{Research partially supported by NSF grants DMS-9800479 (first
author) and DMS-9703775, DMS-0071358 (second author).}

\date{March 27, 2001}

\begin{document}

\begin{abstract}  In their work on  symplectic manifolds, Donaldson and
Auroux use analogues of holomorphic sections of an ample line bundle $L$ over a
symplectic manifold $M$  to create symplectically embedded zero sections and almost
holomorphic maps to
various spaces. Their analogues were termed  `asymptotically holomorphic' sequences
$\{s_N\}$ of sections of
$L^N$. We study another analogue $H^0_J(M, L^N)$ of holomorphic sections, which we call
`almost-holomorphic' sections, following a method introduced earlier by Boutet de Monvel
- Guillemin \cite{BG} in a general setting of symplectic cones. By definition, 
sections in $H^0_J(M, L^N)$  lie in the range of a \szego projector $\Pi_N$.  Starting
almost from scratch, and only using almost-complex geometry, we  construct a simple
parametrix for $\Pi_N$ of precisely the same type as the Boutet de Monvel-Sj\"ostrand
parametrix in the holomorphic case \cite{BS}. We then   show that $\Pi_N(x,y)$  has
precisely the same scaling asymptotics as does the holomorphic \szego kernel as analyzed
in \cite{BSZ1}. The scaling asymptotics imply more or less immediately a number of
analogues of well-known results in the holomorphic case, e.g. a Kodaira embedding
theorem and a Tian almost-isometry theorem.  We also explain how to modify Donaldson's
constructions to prove existence of quantitatively transverse sections in $H^0_J(M,
L^N)$.
\end{abstract}

\maketitle

\tableofcontents

\section*{Introduction}
This paper is concerned with asymptotically
holomorphic sections of ample line bundles
over almost-complex symplectic  manifolds $(M, J, \omega).$  Such line
bundles and sections are symplectic
analogues of the usual objects in complex algebraic geometry.  Interest in
their properties has grown in
recent years because of their use by Donaldson \cite{DON.1, DON.2}, Auroux
\cite{A,A.2} and others \cite{A.3, Sik}
in proving symplectic analogues of standard results in complex geometry.
These results involve properties
of asymptotically holomorphic sections of high powers of the bundle,
particularly those involving their zero sets and the maps they define to
projective space.

We take up the study of asymptotically holomorphic sections from the
viewpoint of the microlocal analysis of
the $\bar{\partial}$ operator on a symplectic almost-complex manifold, and
define a class of `almost holomorphic sections' by a method due to  Boutet
de Monvel and  Guillemin \cite{Bou,BG}. Sections $s$ of powers of a complex
line bundle
$L^N \to M$ over $M$ are identified with equivariant functions $\hat s$ on
the associated $S^1$-bundle $X$, and
the $\bar{\partial}$ operator on sections is identified with the
$\bar{\partial}_b$ operator on $X$.
  In the non-integrable
almost-complex symplectic case there are in general no solutions of
$\bar{\partial}_b {s} = 0$.  To define an `almost
holomorphic section' $s$, Boutet de
Monvel and Guillemin
define  a
certain (pseudodifferential) $\bar{D}$-complex over $X$  \cite{BG}
\cite{BS}. The space $H^0_J(M, L^N)$ of almost holomorphic sections is then
 the space of sections corresponding to   solutions of $\bar{D}_0
{s} = 0$.  The operator   $\bar{D}_0$ is not uniquely or
even canonically defined, and it is difficult to explicitly write down
these almost holomorphic sections.  However, our main results show that on small
length scales  (the same ones used in \cite{DON.1}), these almost holomorphic sections
have the same asymptotic behavior as holomorphic sections in the $N \to \infty$ limit.
Furthermore,  as we shall show in a related paper \cite{SZ2},
typical sections have globally similar properties to the asymptotically
holomorphic sections of Donaldson and Auroux. 

The fundamental ingredient of these statements is   the scaling asymptotics of the almost-holomorphic
\szego kernels $\Pi_N(z, w)$, i.e.
the orthogonal projections onto $H^0_J(M, L^N)$.  As our results will show,
 the `peak sections' $z\mapsto \Pi_N(z,w)$  about points $w\in M$ are
(up to normalization) intrinsically defined almost holomorphic sections
with the main properties of Donaldson's exponentially localized sections
 (cf.\ \cite{DON.1}, Proposition 11). 
Analytically, it is more convenient to work with the lift of $\Pi_N$ as an equivariant kernel 
$\Pi_N(x, y)$ on $X \times X$
than as a section of $L^N \otimes L^{*N}$, so henceforth we will use the term \szego kernel for the equivariant
kernel.  

The scaling asymptotics of $\Pi_N(x,y)$ must be given in special `Heisenberg'
coordinates on $X $, defined as follows: 
 Choose `normal' local coordinates $\{z_j\}$ centered at a point $P_0\in M$ as in
Section
\ref{s-acsymp} and choose a `preferred' local frame for $L$, which together with
the coordinates on $M$ give us `Heisenberg coordinates' $x = (z, \theta), y = (w,
\phi)$ on $X$ (see\/ {\rm \S\ref{s-heisenberg}}).  We obtain the
following near-diagonal asymptotic formula for the \szego kernel in terms of these
Heisenberg coordinates:

\begin{maintheo} \label{neardiag0}  We have:
$$\begin{array}{l}
N^{-m}\Pi_N(P_0+\frac{u}{\sqrtn},\frac{\theta}{N};
P_0+\frac{v}{\sqrtn},\frac{\phi}{N})\\[8pt]  \quad = \frac{1}{\pi^m}
e^{i(\theta-\phi)+u \cdot\bar{v} - \half(|u|^2 + |v|^2)}\left[1+ \sum_{r = 1}^{K}
N^{-\frac{r}{2}} b_{r}(P_0,u,v) + N^{-\frac{K +1}{2}}
R_K(P_0,u,v,N)\right],\\[8pt]  \quad\quad
\mbox{\it where}\ \ 
\|R_K(P_0,u,v,N)\|_{\ccal^j(\{|u|+|v|\le \rho\})}\le C_{K,j,\rho} \ \ \mbox{\it
for}\
\  j=1,2,3,\dots.\end{array}$$ \end{maintheo}

\noindent A more precise statement will be given in Theorem \ref{neardiag}.
The proof of Theorem \ref{neardiag0} is based on  the
construction of explicit parametrices for $\Pi$ and its Fourier coefficients
$\Pi_N$. These parametrices closely resemble those of Boutet de Monvel -
Sj\"ostrand \cite{BS} in the holomorphic case.  The construction is new but
closely follows the work of Menikoff and Sj\"ostrand \cite{MS,Sj} and of
Boutet de Monvel and Guillemin \cite{Bou,BG}.  For the sake of concreteness,
we will give a fairly detailed exposition of the construction of the zeroth
term of the $\bar{D}_j$ complex and of the Szeg\"o kernel.

As more or less immediate corollaries of the scaling asymptotics, we prove the
following symplectic analogues of the holomorphic Kodaira embedding theorem and Tian
almost-isometry theorem
\cite{Ti}. These concern the analogue of the  usual Kodaira map $\Phi_N : M \to
\CP^{d_N-1}$ determined by the space $H^0_J(M, L^N)$. 

\begin{maintheo}\label{tyz} Let  $L \to (M, \omega)$ be the pre-quantum line
bundle
over a $2m$-dimensional
symplectic manifold, and let $\{\Phi_N\}$ be its Kodaira maps.  Then:

\noindent{\rm (a)} There exists a complete asymptotic expansion:
$$ \Pi_N(z,0;z,0)  =  a_0 N^m +
a_1(z) N^{m-1} + a_2(z) N^{m-2} + \dots$$
for certain smooth coefficients $a_j(z)$ with $a_0 = \pi^{-m}$.
 Hence, the  maps $\Phi_N$  are well-defined for $N\gg 0$.

\noindent{\rm (b)}
 Let $\omega_{FS}$ denote the Fubini-Study form on $\CP^{d_N-1}$.
Then $$\|\frac{1}{N}  \Phi_N^*(\omega_{FS}) -
\omega\|_{\ccal^k} = O(\frac{1}{N})$$ for any $k$.  \end{maintheo}

\begin{maintheo}\label{kodaira} The Kodaira map  $\Phi_N:M\to \CP^{d_N-1}$ is an
embedding, for $N$ sufficiently large.
\end{maintheo}

We should draw attention to  related work of   Borthwick-Uribe \cite{BU.1, BU.2}, which defines
`almost holomorphic' in a somewhat different  way, and proves analogues of the Tian almost-isometry
and Kodaira embedding theorems in their setting.   
The  Borthwick-Uribe proof of the almost-complex Tian
theorem was in turn motivated by a
similar proof in the
holomorphic case in \cite{Ze}. The  approach we take here,  based on the off-diagonal scaling asymptotics of the
\szego kernel, is new and we believe it will have further geometric applications.

Finally, we use the 
scaling asymptotics of  the \szego kernel $\Pi_N(x,y)$  to show that
Donaldson's symplectic embedding theorem and the genus calculation of embedded
symplectic curves can be adapted to the setting of almost holomorphic sections.   We
will see in \S 5 that the proof in  \cite{DON.1} of the existence of quantitatively
transversal asymptotically holomorphic sections can be adapted to prove the
existence  of quantitatively transverse sections in $H^0_J(M,L^N)$. 

Further applications of the scaling asymptotics in Theorems \ref{neardiag0} and 
\ref{neardiag} will be given in a forthcoming paper \cite{SZ2}, in which we
discuss  typical properties of almost holomorphic sections (in the probabilistic
sense of \cite{BSZ1, BSZ2}). In particular, 
 we show  that almost every sequence  of
$\lcal^2$-normalized sections $\{s_N\}\in SH^0_J(M,
L^N)$  satisfies the  estimates:
$$\begin{array}{ll} \|s_N\|_\infty = O(\sqrt{\log N}), &\quad \|\nabla^k s_N\|_\infty =
O(N^{k/2}
\sqrt{\log N}),\\ \\  \|\bar{\partial} s_N\|_\infty = O(\sqrt{\log N}), &\quad
\|\nabla^k
\bar{\partial} s_N\|_\infty = O(N^{k/2}\sqrt{\log N}),\quad (k\ge 1)\,.
\end{array}
$$ Consequently, the sequence  of almost-holomorphic sections $\{\frac{1}{\sqrt{\log
N}} s_N\}$ is almost surely asymptotically holomorphic in the sense of Donaldson and
Auroux \cite{DON.1,A}.

  We also use  Theorem \ref{neardiag} to show, in joint
work with P. Bleher
\cite{BSZ2}, that our prior result
\cite[Theorem~3.6]{BSZ1} on the universality of scaling limits of correlations of
zeros of sections of positive line bundles over complex manifolds extends to almost
holomorphic sections on symplectic manifolds.  Thus, on small balls (of radius $1/\sqrt{N}$), almost holomorphic sections have the same  asymptotic probabilities as holomorphic sections of taking on specified values and derivatives.

For the readers' convenience, we conclude this introduction with a brief outline of the
paper. We begin in \S \ref{s-aCR} by describing some terminology from
symplectic geometry and introducing our Heisenberg coordinates. In \S \ref{s-ah}, we
use the method of Boutet de Monvel, Guillemin, and Sj\"ostrand \cite{Bou,
BG, BS} to construct
\szego projectors
$\Pi_N$, which we show are complex Fourier integral operators of the same type
as in the holomorphic case, and then to provide a Boutet de Monvel-Guillemin 
complex of pseudodifferential operators, which replaces the $\dbar_b$ complex in the
symplectic setting. The zeroth term of this complex is used to define sequences of
almost holomorphic sections analogous to the integrable complex case (\S \ref
{s-almostholo}). In \S \ref{s-neardiag}, we  obtain the scaling
asymptotics of $\Pi_N(P_0+\frac{u}{\sqrtn},\frac{\theta}{N};
P_0+\frac{v}{\sqrtn},\frac{\phi}{N})$, which we apply in \S \ref{s-kodaira} to prove
Theorems  \ref{tyz} and \ref{kodaira} on the symplectic Kodaira map. Finally, in \S \ref
{s-transversal}, we show how concentrated
 almost holomorphic sections can be used in place of asymptotically
holomorphic sections to construct quantitatively transversal sections by
Donaldson's method. 

We would like to thank  Michael Christ for informative comments concerning
the off-diagonal decay of the \szego kernel in the holomorphic case.

\section{Circle bundles and almost CR geometry}\label{s-aCR}

We denote by $(M, \omega)$ a compact symplectic manifold such that
$[\frac{1}{\pi}\omega]$ is an integral
cohomology class. 
As is well known (cf.\ \cite[Prop.~8.3.1]{W}; see also \cite{GS}),
there exists a
hermitian line bundle $(L, h) \to M$ and a  metric connection $\nabla$
on
$L$ whose curvature $\Theta_L$ satisfies $\frac{i}{2}\Theta_L = \omega$.
We denote by $L^N$ the
$N^{\rm th}$ tensor power of $L$. The
`quantization' of $(M, \omega)$ at Planck constant $1/N$ should be a Hilbert
space of
polarized sections of $L^N$ (\cite[p. 266]{GS}). In the complex case,
polarized sections are
simply holomorphic sections.
The notion of polarized sections is problematic  in the
non-complex
symplectic setting, since the Lagrangean subbundle $T^{1,0}M$ defining the
complex
polarization is not integrable and there usually are no `holomorphic'
sections.  A subtle but compelling
replacement for the notion of polarized section has been proposed by Boutet
de Monvel and Guillemin \cite{Bou,BG}, and
it is this notion which we adopt in this paper.  For the asymptotic
analysis, it is best to view sections of
$L^N$ as functions on the unit circle bundle $X\subset L^*$; we shall
describe the `almost CR geometry' of  $X$ in \S \ref{s-heisenberg} below.

\subsection{Almost complex symplectic manifolds}\label{s-acsymp}

We begin by reviewing some terminology from almost complex symplectic
geometry. An almost complex
symplectic manifold is a symplectic manifold
$(M,
\omega)$ together with an almost complex structure $J$
satisfying the compatibility condition
$\om(Jv,Jw)=\om(v,w)$ and the positivity condition.
$\omega(v, Jv) > 0$. We give $M$  the Riemannian metric
$g(v,w)=\om(v,Jw)$. We   denote by
$T^{1,0}M,
$ resp.\
$T^{0, 1}M$,  the holomorphic, resp.\ anti-holomorphic, sub-bundle of the
complex tangent bundle $TM$;  i.e., $J = i$ on $T^{1,0}M$ and $J = -i$ on
$T^{0,1}M$.  We give $M$
local coordinates
$(x_1,y_1,\dots,x_m,y_m)$, and we write
$z_j=x_j+iy_j$. As in the integrable (i.e., holomorphic) case, we let
$\{\frac{\d}{\d z_j},
\frac{\d}{\d\bar z_j}\}$ denote the dual frame to $\{dz_j, d\bar z_j\}$.
Although in our case, the coordinates
$z_j$ are not holomorphic and consequently $\frac{\d}{\d z_j}$
is generally not in $T^{1,0}M$, we nonetheless have
$$\frac{\d}{\d z_j} =
\half \frac{\d}{\d x_j}-\frac{i}{2}\frac{\d}{\d y_j}\,,\quad
\frac{\d}{\d\bar
z_j} = \half \frac{\d}{\d x_j}+\frac{i}{2}\frac{\d}{\d y_j}\,.$$
At any point
$P_0\in M$, we can choose  a local frame
$\{\bar Z_1^M,\dots, \bar Z_m^M\}$ for
$T^{0,1}M$ near $P_0$ and coordinates about $P_0$ so that
\begin{equation} \bar Z_j^M= \frac{\d}{\d\bar z_j} + \sum_{k=1}^m B_{jk}(z)
\frac{\d}{\d z_k}\,,\quad B_{jk}(P_0)=0\,,\label{localtangentframe}
\end{equation}  and hence $\d/\d z_j|_{P_0}\in T^{1,0}(M)$. 
This is one of the properties of our `preferred coordinates' defined below.

\begin{defin} Let $P_0\in M$.  A coordinate system $(z_1,\dots,z_m)$ on a
neighborhood $U$ of $P_0$ is {\it preferred\/} at $P_0$ if
$$\sum_{j=1}^m d z_j\otimes d\bar z_j =(g-i\om)|_{P_0} \,.$$
\end{defin}

In fact, the coordinates $(z_1,\dots,z_m)$ are
preferred at $P_0$ if an only if any two of the following conditions (and
hence all three) are satisfied:

\begin{enumerate}
\item[i)] $\quad\d/\d z_j|_{P_0}\in T^{1,0}(M)$, for $1\le j\le m$,
\item[ii)] $\quad\om({P_0})=\om_0$,
\item[iii)] $\quad g({P_0} )= g_0$,
\end{enumerate}
where $\omega_0$ is the standard
symplectic form and $g_0$ is the Euclidean metric:
$$\om_0=\frac{i}{2}\sum_{j=1}^m dz_j\wedge d\bar z_j =\sum_{j=1}^m
(dx_j\otimes dy_j - dy_j\otimes dx_j)\,,\quad g_0=
\sum_{j=1}^m
(dx_j\otimes dx_j + dy_j\otimes dy_j)\,.$$
(To verify this statement, note that condition (i) is equivalent to
$J(dx_j)=-dy_j$ at
$P_0$, and use $g(v,w)=\om(v,Jw)$.)
Note that by the Darboux theorem, we can choose the coordinates so that
condition (ii) is satisfied on a neighborhood of $P_0$, but this is not
necessary for our scaling results.

\subsection{The circle bundle and Heisenberg coordinates}\label{s-heisenberg}

We now let $(M, \omega, J)$ be a compact, almost complex symplectic
manifold such that
$[\frac{1}{\pi}\omega]$ is an integral
cohomology class, and we choose a
hermitian line bundle $(L, h) \to M$ and a  metric connection $\nabla$
on
$L$ with  $\frac{i}{2} \Theta_L = \omega$.
In order to simultaneously analyze sections of all positive powers $L^N$
of the line bundle $L$, we work on the associated principal
$S^1$
bundle $X \to M$, which is defined as follows: let $\pi: L^* \to M$ denote
the
dual line bundle to $L$ with dual metric $h^*$, and put $X = \{v \in L^*:
\|v\|_{h^*} =1\}$.   We let
$\alpha$ be the the connection 1-form on $X$ given by $\nabla$; we then have
$d\alpha =\pi^* \omega$, and thus $\al$ is a contact form on
$X$, i.e., $\al\wedge (d\al)^m$ is a volume form on $X$.

We let $r_{\theta}x =e^{i\theta} x$ ($x\in X$) denote the
$S^1$
action on $X$ and denote its infinitesimal generator by
$\frac{\partial}{\partial \theta}$.
A section $s$ of $L$ determines an equivariant
function
$\hat{s}$ on $L^*$ by the rule $\hat{s}(\lambda) = \left(\lambda,
s(z)
\right)$ ($\lambda \in L^*_z, z \in M$). It is clear that if $\tau \in
\C$
then $\hat{s}(z, \tau \lambda) = \tau \hat{s}$. We henceforth restrict
$\hat{s}$ to $X$ and then the equivariance property takes the form
$\hat{s}(r_{\theta} x) = e^{i \theta}\hat{s}(x)$.  Similarly, a section
$s_N$
of $L^{N}$ determines an equivariant function $\hat{s}_N$ on $X$: put
\begin{equation} \label{sNhat}\hat{s}_N(\lambda) = \left( \lambda^{\otimes
N}, s_N(z)
\right)\,,\quad
\la\in X_z\,,\end{equation} where $\lambda^{\otimes N} = \lambda \otimes
\cdots\otimes
\lambda$;
then $\hat s_N(r_\theta x) = e^{iN\theta} \hat s_N(x)$.
We denote by $\lcal^2_N(X)$ the
space of
such equivariant functions transforming by the $N^{\rm th}$ character.

In the complex case, $X$ is a CR manifold.  In the general almost-complex
symplectic case
it is an almost CR manifold.  The {\it almost CR structure\/} is defined as
follows:
 The kernel of $\alpha$ defines a horizontal hyperplane bundle $H \subset
TX$. Using the projection $\pi: X \to M$, we may lift the splitting
$TM=T^{1,0}M
\oplus T^{0,1}M$ to a splitting $H=H^{1,0} \oplus H^{0,1}$. The almost CR
structure on $X$ is defined to be the splitting
   $TX = H^{1,0} \oplus H^{0,1} \oplus \C \frac{\partial}{\partial
\theta}$.   We  also consider a local orthonormal
frame $Z_1, \dots, Z_n$ of $H^{1,0}$ , resp.\ 
$\bar{Z}_1,
\dots,
\bar{Z}_m$ of $H^{0,1}$, and dual orthonormal coframes $\vartheta_1, \dots,
\vartheta_m,$ resp. $\bar{\vartheta}_1, \dots, \bar{\vartheta}_m$. On the
manifold $X$ we have
$d=
\d_b +\dbar_b +\frac{\partial}{\partial \theta}\otimes \alpha$, where
$\partial_b  =
\sum_{j = 1}^m {\vartheta}_j
\otimes{Z}_j$ and
 $\dbar_b  = \sum_{j = 1}^m \bar{\vartheta}_j \otimes \bar{Z}_j$.
We  define the almost-CR $\bar{\partial}_b$
operator by  $\bar{\partial}_b = df|_{H^{1,0}}$.  
Note that for an $\lcal^2$ section $s^N$ of $L^N$, we have
\begin{equation}\label{dhorizontal}
(\nabla_{L^N}s^N)\nhat = d^h\hat s^N\,,\end{equation} where
$d^h=\d_b+\dbar_b$ is the horizontal derivative on $X$.

Our near-diagonal asymptotics of the \szego kernel (\S \ref{s-neardiag})
are given in terms of the Heisenberg dilations, using local
`Heisenberg coordinates' at a point $x_0\in X$.  To describe these
coordinates, we first need the concept of a
`preferred frame':

\begin{defin}  A {\it preferred frame\/} for $L\to M$ at a point $P_0\in
M$ is a local frame $e_L$ in a neighborhood of $P_0$ such that 

\begin{enumerate}
\item[i)] $\quad \|e_L\|_{P_0} =1$;
\item[ii)] $\quad \nabla e_L|_{P_0} = 0$;
\item[iii)] $\quad \nabla^2 e_L|_{P_0} = -(g+i\om)\otimes
e_L|_{P_0}\in T^*_M\otimes T^*_M\otimes L$.
\end{enumerate} \end{defin}
(A preferred frame can be constructed by multiplying an arbitrary frame by
a function with specified 2-jet at $P_0$; any two such frames agree to
third order at $P_0$.)  Once we have property (ii), property (iii) is
independent of the choice of connection on
$T^*_M$ used to define $\nabla:\ccal^\infty(M,L\otimes T^*_M)\to
\ccal^\infty(M,L\otimes T^*_M \otimes T^*_M)$. In fact, property (iii) is a
necessary condition for obtaining universal scaling asymptotics, because of
the `parabolic' scaling in the Heisenberg group. 
 Note that if $e_L$ is a
preferred frame at $P_0$ and if $(z_1,\dots,z_m)$ are preferred
coordinates at $P_0$, then we compute the Hessian of $\|e_L\|$:
$$\left(\nabla^2
\|e_L\|_h\right)_{P_0} = \Re \left(\nabla^2 e_L,e_L\right)_{P_0}=
-g(P_0)\,;$$ thus if the preferred coordinates are `centered' at
$P_0$ (i.e., $P_0=0$), we have
\begin{equation} \label{asquared}\|e_L\|_h= 1 - \half |z|^2 +
O(|z|^3)\,.\end{equation}

\begin{rem}
Recall (\cite[\S 1.3.2]{BSZ1}) that the Bargmann-Fock representation of the
Heisenberg group acts on the space of holomorphic functions on
$(M,\om)=(\C^m,\om_0)$ that are square integrable with respect to the
weight
$h=e^{-|z^2|}$. We let $L=\C^m\times \C$ be the trivial bundle. Then the
trivializing section
$e_L(z):=(z,1)$ is a preferred frame at $P_0=0$ with respect to the
Hermitian connection $\nabla$ given by $$\nabla e_L=\d \log h \otimes e_L
= -\sum_{j=1}^m \bar z_j dz_j \otimes e_L\,.$$
Indeed, the above yields $\nabla^2e_L|_{0}= -\sum d\bar z_j\otimes dz_j
\otimes e_L(0)= -(g_0+i\om_0) \otimes e_L(0)$.\end{rem}

\medskip
The preferred frame and preferred coordinates together give us `Heisenberg
coordinates':

\begin{defin} A  {\it Heisenberg coordinate chart\/} at a point $x_0$ in
the principal bundle $X$ is a coordinate chart
$\rho:U\approx
V$ with $0\in U\subset \C^m\times \R$ and $\rho(0)=x_0\in V\subset X$ of the
form
\begin{equation}\rho(z_1,\dots,z_m,\theta)= e^{i\theta} a(z)^{-\half}
e^*_L(z)\,,\label{coordinates}\end{equation} where
$e_L$ is a preferred local frame for $L\to M$ at $P_0=\pi(x_0)$, and
$(z_1,\dots,z_m)$ are preferred
coordinates centered at $P_0$.
(Note that $P_0$ has coordinates $(0,\dots,0)$ and $e_L^*(P_0)=x_0$.)
\end{defin}

We now give some computations using local coordinates
$(z_1,\dots,z_m,\theta)$ of the form (\ref{coordinates}) for a local
frame $e_L$.  (For the moment, we do not assume they are Heisenberg
coordinates.) We write \begin{eqnarray*}a(z) &=& \|e^*_L(z)\|^2_{h^*}\ =\ 
\|e_L(z)\|^{-2}_h\,,\\
\alpha &=& d\theta + \beta\,,\qquad \be=\sum_{j=1}^m(A_jdz_j+\bar A_jd\bar
z_j)\,,\\
\nabla e_L&=& \phi \otimes e_L \,,\qquad \mbox{hence}\quad
\nabla e_L^{\otimes N}\ =\ N \phi \otimes e_L^{\otimes N}\,.
\end{eqnarray*} 

We let $\frac{\d^h}{\d z_j}\in H^{1,0}X$ denote the horizontal lift of
$\frac{\d}{\d z_j}$.
The condition $\left(\frac{\d^h}{\d z_j},\al\right) =0$
yields \begin{equation}\label{dhdzj} \frac{\d^h}{\d z_j} = \frac{\d}{\d
z_j} -A_j\frac{\d}{\d
\theta}\,,\quad \frac{\d^h}{\d\bar z_j} = \frac{\d}{\d\bar z_j}
-\bar A_j\frac{\d}{\d\theta}\,.\end{equation}  Suppose $s_N=fe_L^{\otimes
N}$ is a local section of $L^N$. Then by (\ref{sNhat}) and
(\ref{coordinates}),
\begin{equation}\label{sNhat*}\hat s_N(z,\theta) =
f(z)a(z)^{-\half}e^{iN\theta}\,.\end{equation}Differentiating (\ref{sNhat*})
and using (\ref{dhorizontal}), we conclude that
\begin{eqnarray}\phi &=& -\half d\log a -i
\beta\nonumber\\ &=& -\sum_{j=1}^m \left(\half\frac{\d\log a}{\d
z_j}+iA_j\right)dz_j
 -\sum_{j=1}^m \left(\half\frac{\d\log a}{\d \bar
z_j}+i\bar A_j\right)d\bar z_j\,.\label{connection}\end{eqnarray}

Now suppose that $(z_1,\dots,z_m,\theta)$ are Heisenberg coordinates at
$P_0$; i.e., $e_L$ is a preferred frame at $P_0$ and $(z_1,\dots,z_m)$
are preferred coordinates centered at  $P_0$ (with $P_0=0$). By
property (ii) of preferred frames, we have $\phi(0)=0$, and hence by
(\ref{connection})
\begin{equation}\label{da} da|_{0}=d\log a|_{0} =0,\end{equation}
\begin{equation}\label{Aj} A_j(0)=0\,,\quad (1\le j\le
m)\,.\end{equation}
By differentiating (\ref{connection}) and applying the
properties of preferred coordinates and frames, we further obtain
$$\sum_{j=1}^md\bar z_j\otimes dz_j = -\nabla \phi = 
\sum_{j=1}^m d \left(\half\frac{\d\log a}{\d
z_j}+iA_j\right)\otimes dz_j
+\sum_{j=1}^m d \left(\half\frac{\d\log a}{\d \bar
z_j}+i\bar A_j\right)\otimes d\bar z_j\ \mbox{at}\ 0.$$
Thus the following four equations are satisfied at $P_0=0$:
\begin{equation}\label{4equations}\begin{array}{rclrcl}\di
\half\frac{\d^2\log a}{\d z_j \d z_k} +i\frac{\d A_j}{\d z_k} &= & 0\,,
& \di\quad \half\frac{\d^2\log a}{\d z_j \d\bar z_k} +i\frac{\d A_j}{\d\bar
z_k} &= &\de^j_k\,,\\[12pt] \di
\half\frac{\d^2\log a}{\d\bar z_j \d z_k} +i\frac{\d\bar A_j}{\d z_k} &= &
0\,,&\di \quad \half\frac{\d^2\log a}{\d\bar z_j \d\bar z_k} +i\frac{\d\bar
A_j}{\d\bar z_k} &= &0\,,
\end{array}\end{equation} at $P_0$.  Solving (\ref{4equations}) and
recalling that $a(0)=1,\ da|_{0}=0$, we obtain
\begin{equation}\label{d2a}\frac{\d^2 a}{\d z_j \d z_k}(0)=0\,,\qquad
\frac{\d^2 a}{\d z_j \d\bar z_k}(0)=\de^j_k\,,\end{equation}
\begin{equation}\frac{\d A_j}{\d z_k}(0)=0\,,\qquad \frac{\d A_j}{\d\bar
z_k} = -\frac{i}{2}\de^j_k\,.\end{equation}
Hence $A_j=-\frac{i}{2}\bar z_j+O(|z|^2)$ and
\begin{equation}\frac{\d^h}{\d z_j} = \frac{\d}{\d z_j}
+\left[\frac{i}{2}\bar z_j +O(|z|^2)\right]
\frac{\d}{\d \theta}\,,\quad \frac{\d^h}{\d\bar z_j} = \frac{\d}{\d\bar z_j}
-\left[\frac{i}{2} z_j +O(|z|^2)\right]
\frac{\d}{\d \theta}
\,.\label{dhoriz}\end{equation}

\section{Almost holomorphic sections and \szego kernels}\label{s-ah}

In this section, we define the space $H^0_J(M, L^N)$ of almost holomorphic
sections of $L^N$, or equivalently the space ${\mathcal H}^2_N(X)$ of almost CR
functions in
$\hcal^2(X)$ which transform by $e^{i N \theta}$ under the $S^1$ action. The
\szego kernels
$\Pi_N$  are the orthogonal projections onto the spaces ${\mathcal H}^2_N(X)$.

The order of events in the definition of
$H^0_J(M, L^N)$ is almost opposite to that  in the 
holomorphic case.  There, one defines
${\mathcal H}^2(X) $ as the kernel of the 
$\bar{\partial}_b$ operator.  As mentioned above, the kernel is generically empty
in the non-integrable almost complex case.  It is correct but somewhat  misleading
to say that one defines {\it almost holomorphic\/} by deforming 
$\bar{\partial}_b$ in the non-integrable case to a pseudodifferential
operator $\bar D_0$ (which is the zeroth step of a $\bar{D}$-complex), setting
${\mathcal H}^2(X)$ equal to its kernel, and defining
$\Pi$  as the associated orthogonal projection with respect to a metric on $L$.  The misleading aspect is that
actually one  {\it first\/} defines $\Pi$ as a complex Fourier integral projection
operator associated to a certain canonical relation $C$ which naturally generalizes
the one in the integrable case.  Having defined $\Pi,$ one {\it then} defines
$\bar{D}$ as a deformation of $\bar{\partial}_b$ that annihilates $\Pi.$

Existence of a projection $\Pi$ with a prescribed
principal symbol in the algebra $I^*(X \times X, C)$ of Fourier integral kernels
with wave front along $C$ 
 follows from   symbolic and functional analytic arguments, but uniqueness does not.  Indeed, 
$\Pi$ is clearly not unique, since it  could be replaced  by $e^{-i A} \Pi e^{i A}$
for any pseudodifferential operator
$A$ of order $-1$ on $X$. We use the term `\szego kernels' rather than `the
\szego kernel' to remind the reader that
$\Pi$ is by no means  uniquely defined as in the holomorphic case. The
space
$H^0_J(M, L^N)$ therefore also involves  non-unique choices. 

 The non-geometric construction of the spaces of almost holomorphic  sections
may seem less odd if one compares it to the construction of asymptotically
holomorphic sections in
\cite{DON.1}.  There, only
$\bar{\partial}_b$ is defined and asymptotically holomorphic sections are sequences
of sections $\{s_N\}$ which are annihilated by $\bar{\partial}_b$ modulo 
small errors.  For a fixed power $N$, no space of sections is
thereby singled out. On the other hand, our construction gives linear spaces
$H^0_J(M, L^N)$ of almost holomorphic sections, but these spaces depend on the
choice of the projection operator
$\Pi$.

Let us briefly review the holomorphic case from our perspective.  It was
proved in
\cite{BS} that the
\szego kernel of a strictly pseudoconvex domain (including disk bundles of
positively curved line bundles) is a complex Fourier integral operator of the form
$$\Pi(x,y) \sim \int_{\R^+} e^{i t \psi(x,y)} s(x, y ,t) dt,$$ with $s \sim \sum_{n
= 0}^{\infty} t^{m - n} s_n(x,y)$ and with $t \psi(x,y)$ a phase of positive type.
(Complex Fourier integral operators are discussed below.) Note
that the phase is linear in the one-dimensional phase variable
$t$. This is a `hard' result because it shows that an object of complex geometry, the
\szego kernel, has a special singularity structure. It requires estimates
of J. J. Kohn---in particular, the closed range property of $\Box_b$ (cf.\
\cite{BS}).

In
\cite{BSZ1}, we related $\psi$ to the hermitian metric on $L$, which defines the
circle bundle $X$. Let us recall the result. Fix a local holomorphic section $e_L$
of $L$ over $U \subset M$ and define $a
\in \ccal^\infty(U)$ by   $a = |e_L|^{-2}_h$.
Since $L^*|_U \approx U \times \C$ we can define local coordinates on $L^*$
by
$(z, \lambda) \approx \lambda e_L(z)$.
Then a defining function of $X \subset L^*$ is given by $\rho(z, \lambda) =
1 - |\lambda|^2 a(z)$. Define the
function $a(z,w)$ as the almost analytic extension of $a(z)$, i.e.\ the
solution of $\bar{\partial}_z a = 0 =
\partial_w a, a(z,z) = a(z)$ and put $\psi(x,y) = {i}(1 - \lambda
\bar{\mu} a(z,w)).$  Then $t \psi$ is a phase for
$\Pi.$
Because $\Pi$ is a projection, the 
complex canonical relation $C \subset T^*(\tilde{X} \times \tilde{X})$ parametrized by
$t \psi$  must be  idempotent (i.e., 
$C^2 = C = C^*$). As will be explained below, it is the
flow-out of the diagonal in $X \times X$ under the joint Hamiltonian flow of the components of $\bar{\partial}_b$,
which commute by integrability of the complex structure. Such a flow-out is an equivalence
relation ($x \sim y$ if they belong to the same orbit of the commuting flows), and in
this way one also sees that $C$ is idempotent. 

Our purpose now is to construct a complex Fourier integral projection kernel $\Pi$
in the almost complex (non-integrable) case  with the same properties as in the 
complex case. As discussed above, this is not the same kind of `hard' result as in
the complex case because there is no a priori definition of $\Pi$  dictated by almost
complex geometry, so we simply construct one with analogous microlocal properties.
Consequently, we do not need any difficult results of analysis, and the construction
involves only symbol algebra and geometry. Indeed,  the main step is to define the
canonical relation $C$.  Since $\bar{\partial}_b$ is no longer integrable, we
cannot use its components to define the flow-out.  Instead,
we begin in \S \ref{IDEAL} by defining a new positive Lagrangean ideal
$\ical$ corresponding to the ideal generated by the elements of $H^{0,1}X$ in the
complex case.  We then show  in \S \ref{s-C} that $\ical$ defines a unique
canonical relation 
$C$ (which lives in the almost analytic extension of $T^*X\times T^*X$) and we
describe how $C$ can be parametrized by a phase defined on $X\times X\times \R^+$,
as in the holomorphic case.   In \S \ref{s-szego}, we describe the algebra $I^*(X
\times X, C)$ of complex Fourier integral operators with wave front sets along $C$,
and construct a
\szego projector $\Pi\in I^*(X \times X, C)$.  Finally we introduce in \S
\ref{s-almostholo} the Boutet de Monvel - Guillemin $\bar D$ complex of first order
pseudodifferential operators with
$\ker \bar D_0 = \hcal^2(X)$. the operator $\bar D_0$ is the replacement for
$\dbar_b$ in the symplectic case, and the $N^{\rm th}$ Fourier component of
$\hcal^2(X)$ is the space $H^0_J(M,L^N)$ of almost holomorphic sections. 

\subsection{Description of the algebra $I^*(X \times X, C)$}\label{s-I*}

In order to construct an almost holomorphic \szego kernel, we must first
construct the algebra
$I^*(X \times X, C)$ of complex Fourier integral operators in which it lies.  
The definition of
the algebra
$I^*(X \times X, C)$ of complex Fourier integral operators  is essentially given 
in \cite{H} (see Definition
25.5.1), although there it is denoted
$I^*(X \times X, {\mathcal I})$ where ${\mathcal I}$ is a positive conic Lagrangean
ideal.  (We use the older notation of \cite{MeS} and the language of almost analytic
extensions  for ease of comparison with our basic references \cite{MeS, BS, MS}.) It
would take us too far afield to review the definitions and properties of complex
Fourier integral operators, so we only briefly recall some basic ideas and refer the
reader to
\cite{H} or
\cite{MeS} for background.

The most intuitive definition of a Fourier integral kernel $A(x,y) \in I^n(X \times
X, \Lambda)$ is that it is a distribution (generalized function)   defined by a
complex oscillatory integral
$$A(x,y) = \int_{\R^N} e^{i \phi(x, y, \theta)} s(x, y, \theta) d \theta$$ where
$\phi$ is a phase function of positive type parametrizing $\Lambda$, and where $s(x,
y,
\theta)$ is a symbol of order $n+(\dim X
-N)/2$ (cf., \cite{H},
Proposition 25.1.5).  More precisely, we assume $\phi(x,y,
\theta)$ is a regular phase function in the sense of (\cite[Def.~3.5]{MeS}), i.e.
that it has no critical points, is homogeneous of degree one in $\theta$, that the
differentials
$d \frac{\partial \phi}{\theta_j}$ are linearly independent over $\C$ on the
set
$$C_{\phi \R} =\{(x, y, \theta): d_{\theta} \phi = 0\}$$
and such that $\Im \phi \geq 0 $.  
 We then let $\tilde{\phi}(\tilde{x}, \tilde{y}, 
\tilde{\theta})$ be an almost
analytic extension of $\phi$ to $\tilde{X} \times \tilde{X} \times \C^N$.
Here, 
$\tilde{Y}$
denotes the almost analytic extension of a $\ccal^\infty$ manifold $Y$.
 When $Y$ is real
analytic,  $\tilde{Y}$ is the usual complexification of $Y$, i.e. a complex
manifold
 in which $Y$ sits as a totally real submanifold. The reader may find it
simpler to make the extra assumption that $X$ is real analytic. For background on
almost analytic extensions, we
refer to \cite{MeS,MS}.

We 
 put
$$C_{\tilde{\phi}} = \{(\tilde{x}, \tilde{\theta}): d_{\tilde{\theta}}
\tilde{\phi} = 0\}$$
and define the Lagrange immersion
$$\iota_{\tilde{\phi}}: (\tilde{x}, \tilde{\theta}) \in C_{\tilde{\phi}} \to
(\tilde{x},
d_{\tilde{x}} \tilde{\phi}(\tilde{x}, \tilde{\theta})).$$
We say that the  phase $\phi$ parametrizes $\Lambda$ if $\tilde{\Lambda}$ is the image of this
map.

By symbols, we mean the following: the space of classical symbols of order $k$,
denoted 
$S^k(X
\times X
\times
\R^N)$, consists of elements of the form
$s \sim \sum_{j = 0}^{\infty}  s_j(x,y,\theta)$, where 
$s_j$ is a smooth function (defined near $x = y$) of
$(x,y)$ and is homogeneous of order $k-j$ in $\theta\in\R^N$.

Finally, we note that if $\La$ is an equivalence relation, then  $I^*(X \times X,
\La)$ is an algebra \cite{BS,H}.

\subsection{\label{IDEAL} The positive Lagrangean ideal}

We begin with some microlocal analysis of
$\bar{\partial}_b$ in order to introduce the
characteristic variety $\Sigma$
of $\bar{\partial}_b$. In general, we denote by $\sigma_A$ the principal symbol
of a pseudodifferential operator $A$. To describe the principal symbol of
$\bar{\partial}_b$, we   introduce convenient local coordinates and frames.
Recalling that $H X=H^{1,0}X\oplus H^{0,1}X$, we again consider
local orthonormal frames $Z_1, \dots, Z_n$ of $H^{1,0}X$, resp.\
$\bar{Z}_1,
\dots, \bar{Z}_m$ of $H^{0,1}X$, and dual orthonormal coframes $\vartheta_1,
\dots,
\vartheta_m,$ resp. $\bar{\vartheta}_1, \dots, \bar{\vartheta}_m.$ Then we
have $\dbar_b  = \sum_{j = 1}^m \bar{\vartheta}_j \otimes \bar{Z}_j$.  Let
us define complex-valued functions on $T^*X$ by:
$$p_j(x,  \xi) = \langle Z_j(x), \xi),\;\;\;\; \bar{p}_j(x,  \xi) = \langle
\bar{Z}_j(x), \xi \rangle.$$  Then
$$\sigma_{\bar{\partial}_b}(x, \xi) = \sum_{j = 1}^m p_j(x, \xi)
\epsilon(\bar{\vartheta}_j)$$
where $\epsilon$ denotes exterior multiplication.
 We note that $ \{
\bar{p}_j, \bar{p}_k \} = \langle [\bar{Z}_j, \bar{Z}_k],
\xi \rangle$.

To state results, it is convenient to introduce the operator $\boxb :=
\bar{\partial}_b^* \bar{\partial}_b =
 \sum_{j =1}^m \bar{Z}_j^* \bar{Z}_j$ where $\bar{Z}_j^* $ is the
adjoint of the vector field regarded as a linear differential operator.
To conform to the notation of \cite{BG} we also put  $ q=\sigma(\boxb) =
\sum_{j= 1}^m |\bar{p}_j|^2.$
  The
characteristic variety $\Sigma = \{q = 0\}$
of $\bar{\partial}_b$ is the same as that of $\boxb$, namely the vertical
sub-bundle of $T^*X \to M.$
It is the conic submanifold
of $T^*X$  parametrized by  the graph of the contact form,
$\Sigma = \{(x, r \alpha_x): r > 0\} \sim X \times \R^+$.   It
follows that
$\Sigma$ is a symplectic submanifold. It is the dual (real) line bundle
to the vertical subbundle $V \subset TX$, since $\alpha (X) = G(X,
\frac{\partial}{
\partial \theta}).$

We are now ready to introduce  the  positive Lagrangean
ideal $\ical$ whose generators will
define a canonical relation $C$ underlying the \szego kernels and the principal
symbol of the operator $\bar{D}_0$, which will replace $\bar{\partial}_b.$  
 For background on positive Lagrangean
ideals,
see \cite{H}.

\begin{prop}  There exists a unique
positive  Lagrangean ideal $\ical$  with respect to $\Sigma$ containing
$q$. That is, there exists a unique ideal $\ical \subset I_{\Sigma}$ (where
$I_{\Sigma}$ is the ideal of functions vanishing on $\Sigma$) satisfying:
\medskip

\begin{itemize}

\item $\ical$ is closed under Poisson bracket;

\item $\Sigma$ is the set of common zeros
of $f \in \ical$;

\item There exist local generators $\zeta_1, \dots, \zeta_m$ such that the
matrix $\big(\frac{1}{i} \{\zeta_j, \bar{\zeta}_k\}\big)$ is positive
definite
on $\Sigma$ and that $q = \sum_{j,k} \lambda_{j \bar{k}} \zeta_j
\bar{\zeta}_k$, where $\{ \lambda_{j \bar{k}}\}$
is a hermitian positive definite matrix of functions.
\end{itemize}
\end{prop}

\begin{proof}
In the holomorphic case, $\ical$ is generated by the linear
functions
$\zeta_j(x, \xi) =
\langle \xi, \bar{Z}_j \rangle$.  In the general almost complex (or rather
almost CR) setting,
these functions do not Poisson commute and have to be modified. Since the
deviation of
an almost complex structure from being integrable (i.e. a true complex
structure) is measured
by the Nijenhuis bracket, it is not surprising that the generators $\zeta_j$
can be constructed
from the linear functions $\langle \xi, \bar{Z}_j \rangle$ and from the
Nijenhuis tensor. We
now explain how to do this, basically following the method of \cite{BG}.

  As a first approximation to the $\zeta_j$ we begin with the linear
functions $\zeta_j^{(1)}=\bar p_j$ on $T^*X$.  As mentioned above, the
$\zeta_j^{(1)}$ do not
generate a Lagrangean ideal in the non-integrable almost complex case,
indeed
\begin{equation}\label{p} \{\zeta_j^{(1)} , \zeta_k^{(1)} \} =
\langle \xi, [\bar{Z}_j(x), \bar{Z}_k(x)]\rangle\,.\end{equation}
However we do have that
$$\{\zeta_j^{(1)} , \zeta_k^{(1)} \} =
\{\langle \xi, \bar{Z}_j(x)\rangle, \langle \xi, \bar{Z}_k(x)\rangle \} = 0
\;\mbox{on}\; \Sigma.$$ Indeed, for $ (x, \xi) \in \Sigma$, we have $\xi = r
\alpha_x$ for some $r > 0$ so that
\begin{equation}\label{p0}\begin{array}{l}\{\langle \xi,
\bar{Z}_j(x)\rangle, \langle \xi, \bar{Z}_k(x)\rangle\} = r \alpha_x
([\bar{Z}_j(x), \bar{Z}_k(x)]) \\[10pt]\quad= r d \alpha_x (\bar{Z}_j(x),
\bar{Z}_k(x)) = r \pi^* \omega(\bar{Z}_j(x), \bar{Z}_k(x))= 0 \end{array}
\end{equation}
since $\{\bar{Z}_j\}$ forms a Lagrangean subspace for the horizontal
symplectic form $\pi^* \omega$.  Here, $\pi: X \to M$ is the natural
projection.  Moreover if we choose the local horizontal vector fields $Z_j$
to
be orthonormal relative to $\pi^* \omega$, then we also have:
\begin{equation}\label{convex}\begin{array}{l}\{\zeta_j^{(1)} ,
\bar\zeta_k^{(1)} \}(x,\xi) =
\langle \xi, [\bar{Z}_j(x), {Z}_k(x)]\rangle = r \pi^* \omega(\bar{Z}_j(x),
Z_k(x)) \\[10pt] \quad= ir\delta_j^k=i\delta_j^k p_\theta(x,\xi)\,, \qquad
(x,\xi)\in\Sigma\,.\end{array}
\end{equation}
Here, $p_\theta(x,\xi)=\langle \xi,
\frac{\d}{\d\theta}\rangle$.

Finally, we have
$$q = \sum_{j = 1}^m |\langle \xi, {Z}_j\rangle|^2= \sum_{j=1}^m
|\zetaone_j|^2\,.$$ Hence the second and third conditions on the $\zeta_j$
are satisfied by the functions $\zeta^{(1)}_j$.
Furthermore, equation (\ref{p}) tells us that the first condition is
satisfied to zero-th order for the ideal
$\ical_1=\big(\zeta_1^{(1)},\dots,\zeta_m^{(1)}\big)$.  In fact, let us
precisely describe the error.  We consider the orthonormal (relative to
$\om$) vector fields $Z_j^M=\pi_*Z_j$ of type (1,0) on $M$.  Recall that
the Nijenhuis tensor is given by
$$N(V,W)=\half\big([JV,JW]-[V,W]-J[V,JW]-J[JV,W]
\big)\,.$$
Hence,
\begin{equation}\label{Nijen}
N(Z_j^M,Z_k^M)=(-1-iJ)[Z_j^M,Z_k^M]=-2[Z_j^M,Z_k^M]_{(0,1)}
\eqd\sum_{p=1}^m
N^p_{jk}\bar Z_p^M\,.\end{equation}
We note that by definition, \begin{equation}\label{Nijensym}
N^p_{jk}=N^p_{kj}\,.
\end{equation} Furthermore, by the Jacobi identity
$$\{\{\zeta_j,\zeta_k\},\zeta_p\}
+\{\{\zeta_p,\zeta_j\},\zeta_k\} +\{\{\zeta_k,\zeta_p\},\zeta_j\}=0$$
applied to $(x,\al_x)\in\Si$, we have
\begin{equation}\label{Nijensym1}
N^p_{jk}+N^k_{pj}+N^j_{kp}=0\,.
\end{equation}

By (\ref{p0}) and (\ref{Nijen}), we have
\begin{equation}\label{pi} \{\zeta_j^{(1)} , \zeta_k^{(1)} \} =
\sum_{p=1}^m f_p^1\zeta_p^{(1)} +  \sum_{p=1}^m\bar N^p_{jk}
\bar\zeta_p^{(1)}\,.\end{equation}

We now argue, following \cite{BG}, that these functions can be successively
modified to satisfy the same conditions to infinite order on $\Sigma.$
The next step is to modify the functions $\zetaone_j$ by  quadratic terms so
that they satisfy the conditions $\{\zeta_j,\zeta_k\}\in \ical$ to first order
and the condition $q =
\sum_{j} |\zeta_j|^2$ to order 3 on
$\Sigma$. So we try to construct
new functions
$$\zetatwo_p=\zetaone_p + R_p\,,\quad R_p=\sum_{j,k}\nu_p^{jk}\bar\zetaone_j
\bar\zetaone_k$$   so that
\begin{eqnarray}\{\zetatwo_j , \zetatwo_k \} &=& \sum_{p}
f_p^2\zetatwo_p + \sum_{\al_1,\al_2}
\mu_{jk}^{\al_1\al_2}\bar\zetatwo_{\al_1}
\bar\zetatwo_{\al_2}\,;\label {i2}\\
q &=& \sum_{p}v_p^2 \zetatwo_p +\sum_\al \phi_p^\al \bar\zetatwo_{\al_1}
\bar\zetatwo_{\al_2} \bar\zetatwo_{\al_3}\bar\zetatwo_{\al_4}\,,\quad
(\al=(\al_1,\dots,\al_4))
\,.\label{ii2}\end{eqnarray}

Let us now solve (\ref{i2})--(\ref{ii2}) for the $\nu_p^{jk}$.  First of
all,
we choose $\nu_p^{jk}=\nu_p^{kj}$.  We have
$$\{\zetatwo_j,\zetatwo_k\}=\sum_{p=1}^m f_p^1\zeta_p^{(1)} +
\sum_{p=1}^m\bar
N^p_{jk}\bar\zetaone +\{\zetaone_j,R_k\} -\{\zetaone_k,R_j\} \mod
I_\Sigma^2\,.$$
By (\ref{convex}), we have
\begin{equation}\label{convex2} \{\zetaone_j,\bar\zetaone_k\}=i\de^k_j
p_\theta
\mod I_\Si\,, \end{equation} and thus
$$\{\zetaone_j,R_k\}=\sum_{p=1}^m 2i\nu_k^{pj}p_\theta\bar\zetaone_p \mod
I_\Si^2\,.$$ Therefore,
\begin{equation}\label{p2}
\{\zetatwo_j,\zetatwo_k\}= \sum_{p=1}^m f_p^1\zetatwo_p +\sum_{p=1}^m
\left(\bar N_{jk}^p +2i(\nu^{pj}_k-\nu^{pk}_j)p_\theta\right) \bar\zetaone_p
\mod I_\Si^2\,.\end{equation}
Hence $$\bar N_{jk}^p=
2i(\nu^{pk}_j-\nu^{pj}_k)p_\theta \quad \mbox{on}\ \Si\,,$$
or equivalently,
\begin{equation}\label{c2}
\nu^{pk}_j-\nu^{pj}_k =\frac{i}{2p_\theta}\bar N_{jk}^p \mod I_\Sigma\,.
\end{equation}

On the other hand, \begin{eqnarray*}q &=& \sum_p |\zetatwo_p-R_p|^2 \ =\
\sum_pv_p^2\zetatwo_p -R_p\bar\zetatwo_p\\
&=& \sum_pv_p^2\zetatwo_p  -\sum_{j,k,p} \nu_p^{jk} \bar\zetatwo_j
\bar\zetatwo_k
\bar\zetatwo_p + \sum_\al \phi_p^\al \bar\zetatwo_{\al_1}
\bar\zetatwo_{\al_2} \bar\zetatwo_{\al_3}\bar\zetatwo_{\al_4}\,.
\end{eqnarray*}
Hence (\ref{ii2}) is equivalent to
\begin{equation}\label{nusym} \nu_p^{jk}+\nu_k^{pj}+\nu_j^{kp}=0\,.
\end{equation} Using (\ref{Nijensym})--(\ref{Nijensym1}), we can solve
the equations (\ref{convex2}) and (\ref{c2}) to obtain
\begin{equation}\label{soln2} \nu_p^{jk}=\frac{i}{6p_\theta} \left(\bar
N^k_{pj} + \bar N^j_{pk}\right)\,.\end{equation}
Indeed, the solution (\ref{soln2}) is unique (modulo $I_\Sigma$) and hence
the
$R_p$ are unique modulo $I_\Sigma^3$.  In summary,
\begin{equation} \zetatwo_p=\zetaone_p + \frac{i}{3p_\theta}\sum_{j,k}
\bar N^k_{pj}\bar\zetaone_j \bar\zetaone_k\,.\label{zetatwo}\end{equation}

The passage from the $n^{\rm th}$ to the $(n+1)^{\rm st}$ step is similar,
and we refer to \cite[pp.~147--149]{BG}. \end{proof}

\begin{rem} Define $p_{\theta}(x, \xi) = \langle \xi,
\frac{\partial}{\partial \theta}\rangle.$ Since
the joint zero set of $\{\zeta_1, \dots, \zeta_m\}$ equals $\Sigma$ and
since $p_{\theta} \not= 0$ on
$\Sigma - 0$ it follows that $\{\zeta_1, \dots, \zeta_m, p_{\theta}\}$ is an
elliptic system of symbols.
\end{rem}

\subsection{The complex canonical relation}\label{s-C}

We now construct the underlying complex canonical relation $C$ and we show
that $C$ can be parametrized by a phase of the form $t\psi(x,y)$ as in the complex
case (Theorem~\ref{oscintth}). 

As mentioned above,
$C$ does not live in $T^*
X \times T^*X$ but rather in
its almost analytic extension $T^* \tilde{X} \times T^* \tilde{X}$. 
Since $\pi: X \to M$ is an $S^1$ bundle over $M$, its complexification
$\tilde{\pi}: \tilde{X} \to \tilde{M}$
defines a $\C^*$ bundle over $\tilde{M}$. The connection form $\alpha$ has
an (almost) analytic continuation to
a connection $\tilde{\alpha}$ to this bundle and we may split $T \tilde{X} =
\tilde{H} \oplus \tilde{V}$,
where $\tilde{V} \to T \tilde{M}$ is the vertical subbundle of the fibration
$\tilde{X} \to \tilde{M}$ and where
$\tilde{H} \to T \tilde{M}$ is the kernel of $\tilde{\alpha}.$

The (almost) complexification of $T^* X$ is of course $T^*(\tilde{X})$. We
denote the canonical symplectic form
on $T^*X$ by $\sigma$ and that on $T^*(\tilde{X}$ by $\tilde{\sigma}$; the
notation is consistent because it is
the complexification of $\sigma.$  The symplectic cone $\Sigma$ complexifies
to $\tilde{\Sigma}$ and it remains
symplectic with respect to $\tilde{\sigma}.$ It is given by $\{(\tilde{x},
\tilde{\lambda} \tilde{\alpha}_{\tilde{x}}):
\tilde{\lambda} \in \C^*\}.$  We have a natural identification $L^* \iff
\Sigma$ given by $r x \to (x, r \alpha_x).$
We further note that the $\C^*$ bundle $L^* \to M$ is the fiberwise
complexification of the $S^1$ bundle $X \to M$,
hence $L^* \to M$ is the restriction of $\tilde{\pi}$ to $\tilde{\pi}^{-1}
(M).$ We will therefore view $L^*$ as
a submanifold of $\tilde{X}.$

We let $\tilde{\zeta}_j$ be almost analytic extensions to $\tilde X$  of the functions
$\zeta_j$.  We then define
\begin{equation} \mathcal{ J}_+ = \{(\tilde{x}, \tilde{\xi}) \in T^*
\tilde{X}: \tilde{\zeta}_j = 0\  \forall j\}\;, \end{equation}
which is an involutive manifold of $T^*\tilde{X}$ with the properties:
\begin{equation} \begin{array}{rl} \mbox{(i)} & (\jcal _+)_{\R} = \Sigma \\ & \\
\mbox{(ii)} & q |_{\jcal _+} \sim 0 \\ &  \\\mbox{(iii)} &
\frac{1}{i} \sigma(u, \bar{u}) > 0, \forall u \in T(\jcal _+)^{\bot}\\ &
\\\mbox{(iv)} &
T_{\rho}(\jcal _+) = T_{\rho} \tilde{\Sigma} \oplus W_{\rho}^+.
\end{array} \end{equation}
Here, $W_{\rho}^+$ is the sum of the eigenspaces of $F_{\rho}$,
the normal Hessian
of $q$,  corresponding to the eigenvalues
$\{ i \lambda_j\}$. The null foliation of $\mathcal{J}_+$ is given by
the joint Hamilton flow of the $\tilde{\zeta}_j$'s.

The following proposition,  proved in \cite{MS} and in
(\cite{BG}), Appendix, Lemma 4.5), adapts to our context and  defines the complex
canonical relation $C$:

\begin{prop} There exists a unique strictly positive almost analytic
canonical relation
$C$
satisfying
$$ \diag(\Sigma) \subset C \subset \jcal _+ \times
\overline{\jcal _+}.$$  \label{C}\end{prop}

Indeed,
\begin{equation} C = \{(\tilde{x}, \tilde{\xi}, \tilde{y}, \tilde{\eta}) \in
\jcal _+ \times \overline{\jcal _+}: (\tilde{x}, \tilde{\xi}) \sim
(\tilde{y}, \tilde{\eta})\}, \end{equation}
where $\sim$ is the equivalence relation of `belonging to the same leaf of
the null foliation of $\jcal _+.$  Thus, $C$ is the flow-out of its real
points, $\diag(\Sigma)$, under the joint Hamilton
flow of the $\tilde{\zeta}_j$'s. It is clear from the description that $C
\circ C = C^* = C,$ i.e. that $C$
is an idempotent canonical relation.

We now prove that  $C$ can
be parametrized by a phase $t\psi(x,y)$
defined
on $ X \times X\times \R^+$.  This is helpful in analyzing the scaling limit
of $\Pi_N(x,y)$.
In the following we use local coordinates $(z, \lambda)$ on $L^*$ coming
from a choice of local coordinates $z$ on $M$ and a local frame $e_L(z)$ of
$L$,
and a corresponding local trivialization $(\tilde{z}, \lambda)$ of
$\tilde{X} \to \tilde{M}$. As before, we let $a=\|e_L^*\|^2$.

\begin{theo} \label{oscintth}  There exists a unique regular phase function $ t
\psi(x,y)
\in \ccal^\infty( X \times X\times \R^+ )$  of positive type 
such that:
\begin{itemize}

\item
$id_x \psi |_{x = y}= -i d_y \psi
|_{x = y} = \alpha$;

\item  The 
 almost analytic extension  $\tilde{\psi} \in \ccal^\infty(\tilde{X}
\times \tilde{X})$ of $\psi$ has the form
$\tilde\psi(\tilde x,\tilde y) =
i( 1 - \lambda \bar{\mu} \tilde a(\tilde z,\tilde w))$
 with $\tilde a(z,z)=a(z)$ and $\tilde a(\tilde z,\tilde w) =
\overline{\tilde a(\tilde w,\tilde z)}$;

\item  $t\tilde{\psi}$ parametrizes $C$.

\end{itemize} 
\end{theo}

\begin{proof}

We need to  construct a function $a(z,w)$ so that $  t \tilde\psi$ as above
parametrizes the
canonical relation $C$, i.e.\ that $C$ is the image of the Lagrange immersion
\begin{equation}\begin{array}{l}  \iota_{\tilde{\psi}}: C_{t\tilde{\psi}}= 
\R^+
\times \{\tilde{\psi} = 0\}\to T^*(\tilde{X} \times \tilde{X})\\ \\
(t, \tilde{x}, \tilde{y}) \mapsto (\tilde{x}, t d_{\tilde{x}} \tilde{\psi};
\tilde{y}, - t
d_{\tilde{y}} \tilde{\psi}) \end{array} \end{equation}
 Since $C$ is the unique canonical relation satisfying
$ \label{C*} \diag(\Sigma) \subset C \subset \jcal _+ \times
\overline{\jcal _+}$,
 the conditions that  $t\tilde{\psi}$ parametrize $C$ are the following:
\begin{enumerate}
\item  [{\rm i)}] $\{(x,y)\in X\times X:{\psi}(x,y) = 0\} = \diag(X)$;
\item  [{\rm ii)}] $d_x \psi |_{x = y}= - d_y \psi
|_{x = y} = r \alpha$ for $x, y \in X$ and for some function $r(x) > 0$;
\item   [{\rm iii)}]   $\tilde{\zeta}_j(\tilde x, d_{\tilde x} \tilde{\psi})=
0 =
\tilde{\zeta}_j(\tilde y, d_{\tilde y} \tilde{\psi})$ on $\{\tilde{\psi} =
0\}.$
\end{enumerate}
Such a $\tilde{\psi}$ is not unique, so we  require that $r\equiv 1$ in
condition (ii), i.e., 
$$d_x \psi |_{x = y}= - d_y \psi
|_{x = y} = \alpha\,.$$

Suppose we have $\tilde\psi(\tilde x,\tilde y) =
i( 1 - \lambda \bar{\mu} \tilde a(\tilde z,\tilde w))$. We observe that
$$\tilde{\psi} = 0 \iff  \tilde{a}(\tilde{z}, \tilde{w}) = (\lambda
\bar{\mu})^{-1},$$ and hence
\begin{equation} \label{first} \begin{array}{l}   id_{\tilde{x}} \tilde{\psi}
=  \bar{\mu}
\tilde{a}(\tilde{z}, \tilde{w}) {d \lambda} +
 \lambda \bar{\mu} d_{\tilde{z}} \pi^* \tilde{a}(\tilde{z}, \tilde{w}) \\ \\
\quad = {\lambda} ^{-1}{d \lambda} +
  \tilde{a}^{-1}  d_{\tilde{z}} \tilde{a} (\tilde{z}, \tilde{w})\;\;
\iff \;\; \tilde{\psi} = 0. \end{array} \end{equation}
The conditions on $a$ are therefore:
$$\left\{ \begin{array}{ll}
 a(z,w)\lambda \bar{\mu} = 1 \iff (z, \lambda ) = (w, \mu) \in X; \\ \\
(a\inv d_z  a + \la\inv d\la)|_{\diag(X)} =- (a\inv d_w  a + \la\inv
d\la)|_{\diag(X)} =  \alpha\\ \\
\tilde{\zeta}_j\left(\tilde{z}, {\lambda},\la\inv d\la +
 \tilde{a}^{-1}  d_{\tilde{z}} \tilde{a} (\tilde{z}, \tilde{w})\right) = 0 =
\tilde{\zeta}_j\left(\tilde{w}, \tilde{\mu}, \mu\inv{d \mu} +
 \tilde{a}^{-1}  d_{\tilde{w}} \tilde{a} (\tilde{z}, \tilde{w})\right),\
\forall (z,w, \lambda, \mu)
 \end{array} \right.$$

A solution $a(z,w)$ satisfying the first condition must satisfy  $a(z,z)
|\lambda|^2 = 1$ on $X$, so that $a(z,z) |\lambda|^2$ is the
local hermitian metric on $L^*$ with
unit bundle $X$, i.e. $a(z,z)=a(z)$.

We now prove that these conditions have  a unique solution near the
diagonal. We do this by
  reducing the canonical relation $C$ by
the natural  $S^1$ symmetry. The reduced relation $C_r$ has a unique
generating
function $\log a$;  the three conditions  above on $a$ will follow
automatically from this fact.

 The $S^1$ action of $X$ lifts to $T^*X$ as the Hamiltonian flow
of the function $p_{\theta}(x, \xi):= \langle \xi, \frac{\partial}{\partial
\theta} \rangle.$ The $\zeta_j$
 are invariant under this $S^1$ action, hence
\begin{equation}\label{PC} \{p_{\theta}, \zeta_j \} = 0\;\;\forall j.
\end{equation}

Now consider the level set $\{p_{\theta} = 1\} \subset T^*X$.  Dual to the
splitting $TX = H \oplus V$ we get
a splitting $T^*X = H^* \oplus V^*$, where
$$V^*(X) = \R \alpha = H^o,\;\;\;\; H^*(X) = V^o$$
where $E^o$ denotes the annihilator of a subspace $E$, i.e. the linear
functionals which vanish on $E$. Thus,
$p_{\theta} = 0$ on the horizontal space $H^*(X)$ and $p_{\theta}(\alpha) =
1$. Since
 $p_{\theta}$ is linear on the fibers of $T^*X$,
the set $\{p_{\theta} = 1\}$ has the form $\{\alpha + h: h \in H^*(X)\}$.
We also note that $p_{\theta} (d \theta) = 1$ in the  local coordinates
$(z, \theta)$ on $X$ defined by  $\lambda = e^{i \theta}$.   Hence
$\{p_{\theta} = 1\}$  may also be
identified with  $\{d \theta + h: h \in H^*(X)\}$.

Since $\{p_{\theta} = 1\}$ is a hypersurface, its null-foliation is given by
the orbits of the Hamiltonian
flow of $p_{\theta}$, i.e. by the $S^1$ action.  We use the term `reducing
by the $S^1$-action' to mean setting
$p_{\theta} = 1$ and then dividing by this action.  The reduction of $T^*X$
is thus defined by $(T^*X)_r = p_{\theta}^{-1}(1)/S^1$.  
Since $p_{\theta}^{-1}(1)$ is an affine bundle over $X$ with fiber
isomorphic to $H^*(X)\approx T^*M,$ it is clear that $(T^*X)_r\approx T^*M$ as
vector bundles over $M$.  We can obtain
a symplectic equivalence using the local coordinates $(z, \theta)$ on $X$.
Let $(p_z, p_{\theta})$ be the corresponding
symplectically dual coordinates, so that the natural symplectic form
$\sigma_{T^*X}$ on $T^*X$ is given by $\sigma_{T^*X} =
dz \wedge d p_z + d\theta \wedge d p_{\theta}$.  The notation $p_{\theta}$
is consistent with the above. Moreover,
the natural symplectic form on $T^*M$ is given locally by $\sigma_{T^*M} =
dz \wedge d p_z.$
Now define the projection
$$\chi: p_{\theta}^{-1}(1) \to T^*M, \;\;\; \chi(z, p_z, 1, p_{\theta}) =
(z, p_z). $$
This map commutes with the $S^1$ action and hence descends to the quotient
to define a local map  over $U$, still denoted $\chi$, from $(T^*X)_r \to
T^*M$. Clearly $\chi$ is symplectic.

We now reduce the canonical relation $C$. Thus we consider the $\C^* \times
\C^*$ action on $T^*\tilde{X} \times T^*\tilde{X} - 0$ generated by
$p_{\theta}(x, \xi), p_{\theta}(y, \eta).$  The reduction of $C$ is given by
$$C_r = C \cap (p_{\theta} \times p_{\theta})^{-1}(1,1) / \C^* \times
\C^*.$$  We then use $\chi \times \chi$ to  identify $C_r$ with a
(non-homogeneous) positive canonical relation in $T^*(\tilde{M} \times
\tilde{M}).$  Thus in coordinates,
\begin{equation} C_r = \{(\tilde{z}, \tilde{p_z} , \tilde{w},
\tilde{p_w}) \in T^*(\tilde{M} \times \tilde{M}): \exists \lambda, \mu,
(\tilde{z}, \lambda,  \tilde{p_z} , 1;  \tilde{w}, \mu, \tilde{p_w}, 1) \in
C\}.
 \end{equation}

Since reduction preserves real points, it is clear that
$$\begin{array}{l} (C_r)_{\R} = C_{\R} \cap (p_{\theta} \times
p_{\theta})^{-1}(1,1) / \C^* \times
\C^* \\ \\
= \{(z, p_z , z,p_z) \in \diag(T^*(M \times M)): \exists \theta \ \mbox{such
that}\
\alpha_{z, e^{i \theta}} = d\theta + p_z\}. \end{array} $$

Let us denote by $\tilde{\zeta}_{j r}$ the reductions of the functions
$\tilde{\zeta}_j$
by the $S^1$ symmetry.  Then $\tilde{\zeta}_{j r} = 0$ on either pair of
cotangent vectors
in $C_r$.  Moreover, by the uniqueness statement on $C$ it follows that
$C_r$ is the unique
canonical relation in $T^*(\tilde{M} \times \tilde{M})$ with the given set
of real points
and in the zero set of the $\tilde{\zeta}_{j r}$'s.

We now observe that $C_r$ has, at least near the diagonal, a unique global
generating
function. This holds because the natural projection
\begin{equation} C_r \subset T^*(\tilde{M} \times \tilde{M})  \to  \tilde{M}
\times \tilde{M} \end{equation}
is a local diffeomorphism near the diagonal.  Indeed, its derivative gives a
natural isomorphism
\begin{equation} T_{\rho, \rho} C_r \approx H^* \oplus H^* \approx T(\tilde{M}
\times \tilde{M})\,. \end{equation}
Therefore, there exists a  global generating function $\log \tilde{a} \in
\ccal^\infty(\tilde{M} \times \tilde{M})$
i.e.
\begin{equation} C_r =\{(\tilde{z}, d_{\tilde{z}} \log \tilde{a}, \tilde{w},
d_{\tilde{w}} \log \tilde{a}), \;
\tilde{z}, \tilde{w} \in \tilde{M}\} . \end{equation}
Since $C^* = C$ it follows that $C_r^* = C^r$ and hence that $a(w,z) =
\overline{a(z,w)}.$

Working backwards, we find that
 the function $\tilde \psi(\tilde x,\tilde y) = i(1 -
\lambda \bar{\mu} \tilde a(\tilde z,\tilde w))$ satisfies the equations
$\tilde{\zeta}_j(\tilde{x}, d_{\tilde{x}} \tilde{\psi}) =
\tilde{\zeta}_j(\tilde{y}, d_{\tilde{y}} \tilde{\psi}) = 0$ on $\tilde{\psi}
= 0.$ Therefore the Lagrange immersion
\begin{equation}\begin{array}{l}  i_{\tilde{\psi}}: C_{t\tilde{\psi}}=  \R^+
\times \{\tilde{\psi} = 0\}\to T^*(\tilde{X} \times \tilde{X})\\ \\
(t, \tilde{x}, \tilde{y}) \to (\tilde{x}, t d_{\tilde{x}} \tilde{\psi};
\tilde{y}, - t
d_{\tilde{y}} \tilde{\psi}) \end{array} \end{equation}
takes its image inside $\jcal _+ \times \overline{\jcal _+}$ and reduces to
$C_r$ under
the $S^1$-symmetry.  To conclude the proof it is only necessary to show that
the real
points of the image of $i_{\tilde{\psi}}$ equal $\diag(\Sigma).$ We know
however that these
real points reduce to $(C_r)_{\R}$ and hence that $z = w$ at real points.
But we have
$$1 = \lambda \bar{\mu} a(z,w) = e^{i (\theta - \phi)}
\frac{a(z,w)}{\sqrt{a(z)} \sqrt{a(w)}},\;\;\;
\mbox{on}\;\; \{\tilde{\psi} = 0\}$$
hence when $z = w$ we have $e^{i(\theta - \phi)} = 1$ and hence $x = y$.
Since
$d_{\tilde{x}} \tilde{\psi}(x,y)|_{x = y} = \alpha_x$, it follows that the
real points
indeed equal $\diag(\Sigma)$. Therefore $t \tilde{\psi}$ parametrizes $C$.

To show that the phase is of positive type, we need to describe the
asymptotics of $a(z,w)$ near the diagonal.  Note that
in the almost-complex case, we cannot describe
$a(z,w)$ as the almost analytic extension of $a(z,z)$.
(Of course, $\tilde a(\tilde z,\tilde w)$ is the almost
analytic extension of $a(z,w)$, by definition.) For our near-diagonal
asymptotics in the nonintegrable case, we instead use the following second
order expansion of
$a$ at points on the diagonal:

\begin{lem} Suppose that $(z_1,\dots,z_m)$ are preferred coordinates and
$e_L$ is a preferred frame at a point $P_0\in M$. Then the Taylor expansion of
$a(z,w)$ at $z=w=0$ is
$$a(z,w) = 1 + z\cdot \bar w + \cdots\;.$$
\label{a2}\end{lem}

\begin{proof} To begin, we recall that $a(0,0)=a(0)= \|e_L^*(P_0)\|^2=1$. To
compute the first and second order terms, we return to the equation
\begin{equation}\label{forallzw}\zeta_j\big(z,\la, \frac{d \lambda}{\lambda}
+ d_z
\log a(z,w)\big) = 0,
\;\;\;\;
\forall (z,\la;w) \in X \times M.\end{equation}
Let us write $\zeta_j= \zeta_j^{(1)}
+ R_j^{(2)},$ where  $R_j^{(2)}$ vanishes to second order on $\Sigma$ and we
recall that $\zetaone_j(\xi)=(\bar Z_j,\xi)$. Let us also Taylor
expand $\log a$:
$$\log a = L(z,w)+Q(z,w) +\cdots\,,$$ where $L$ is linear and $Q$ is
quadratic. Since $e_L$ is a preferred frame at $P_0$, it follows from
(\ref{asquared}) that
$a(z,z)=1+|z|^2 +\cdots$ and hence 
\begin{equation}\label{logadiag} L(z,z)=0\,,\qquad
Q(z,z)=|z|^2\,.\end{equation}

Since $d_z \log a|_{z = w}  +  \frac{d \lambda}{\lambda} =
\alpha\in\Sigma$, it follows from (\ref{forallzw}) that
\begin{equation}\label{Rj2}\zetaone_j\big(z,\la,\frac{d\lambda}{\lambda} +
d_z \log a\big) = -R_j^{(2)}(z, \lambda, \frac{d \lambda}{\lambda} +
d_z
\log a) = O(|z- w|^2).\end{equation} Since $a(z,w) =\overline{a(w,z)}$, we
can write
$$L(z,w)=\sum_{j=1}^m(b_jz_j+c_j\bar z_j+\bar c_jw_j+\bar b_j\bar w_j)\,.
$$ Since the $z_j$ are preferred coordinates and
$e_L$ is a preferred frame at $P_0$, we can choose the $\bar Z_j$ so that
$\bar Z_j(0)=\frac{\d}{\d\bar z_j}$ and hence by (\ref{Rj2}),
$$0=\zetaone_j\left.\big(z,\la,\frac{d\lambda}{\lambda} + d_z \log
a\big)\right|_{z=w=0,\la=1} =
\left.\left(\frac{\d}{\d
\bar z_j}, d_z \log a\right)\right|_{(0,0)}= c_j\ \forall j\,.$$  Since
$L(z,z)=0$, we have
$b_j+\bar c_j = 0$, and hence $L=0$.

To investigate the quadratic term $Q$ in (\ref{logadiag}), we write
\begin{equation}\label{findQ}(\frac{d\la}{\la}+d_z\log a)|_{(z,w)} 
=\alpha_{z} +\sum_{j=1}^m\left[z_j U'_j+\bar z_j U''_j+ w_jV'_j + \bar
w_jV''_j\right] + O(|z|^2+|w|^2)\,,
\end{equation} where  
$$\begin{array}{lcllcl}U'_j &=& \sum_{k=1}^m\left(\frac{\d^2 Q}{\d z_j
\d z_k} dz_k + \frac{\d^2 Q}{\d z_j \d\bar z_k}d\bar z_k\right)\,,\quad &
U''_j &=& \sum_{k=1}^m \left(\frac{\d^2 Q}{\d\bar z_j
\d z_k} dz_k +  \frac{\d^2 Q}{\d\bar z_j \d\bar
z_k}d\bar z_k\right)\,,\\[10pt]
V'_j &=& \sum_{k=1}^m\left(\frac{\d^2 Q}{\d w_j
\d z_k} dz_k + \frac{\d^2 Q}{\d w_j \d\bar z_k}d\bar z_k\right)\,,\quad &
V''_j &=& \sum_{k=1}^m \left(\frac{\d^2 Q}{\d\bar w_j
\d z_k} dz_k +  \frac{\d^2 Q}{\d\bar w_j \d\bar
z_k}d\bar z_k\right)\,.
\end{array}$$

Applying $\zetaone_k$ to (\ref{findQ}) and using (\ref{Rj2}) and the fact that 
$\zetaone_k(\al_z)=0$, we have 
\begin{equation}\sum_{j=1}^m \left[z_j(\bar Z_k|_z,U'_j) + \bar z_j
(\bar Z_k|_z,U''_j) + w_j(\bar Z_k|_z,V'_j) + \bar w_j
(\bar Z_k|_z,V''_j)\right]=O(|z|^2+|w|^2)\,.\label{barZk}\end{equation}
By (\ref{localtangentframe}) and (\ref{dhdzj}),
$$ \bar Z_k|_z= \frac{\d}{\d\bar z_k} + \sum_{l=1}^m B_{kl}(z)
\frac{\d}{\d z_l} +  C_k(z)\frac{\d}{\d
\theta} \,,\quad B_{kl}(0)=0\,.$$ Hence by (\ref{barZk}),
$$\frac{\d^2 Q}{\d z_j\d\bar z_k}=\left(\frac{\d}{\d\bar z_k},U'_j\right)=
(\bar Z_k|_0,U'_j)=0\,.$$  Similarly, $\frac{\d^2 Q}{\d\bar z_j\d\bar z_k}=
\frac{\d^2 Q}{\d w_j\d\bar z_k} = \frac{\d^2 Q}{\d \bar w_j\d\bar z_k} =0$.
Thus $Q(z,w)$ has no  terms containing $\bar z_k$. 
Since $Q(z,w)=\overline{Q(w,z)}$, the quadratic function
$Q$ also has no terms containing $w_k$, so we can write
$$Q(z,w)=B(z,z)+H(z,\bar w)+\overline{B(w,w)}\,,$$ where $B$, resp.\ $H$, is a
bilinear, resp.\ hermitian, form on $\C^m$. Since $Q(z,z)=|z|^2$ (recall
(\ref{logadiag})), we conclude that $B(z,z)=0$ and
hence $Q(z,w)=H(z,\bar w)=z\cdot\bar w$.
\end{proof}

We now complete the proof that
the phase is of positive type; i.e., $\Im\psi\ge 0$ 
on some neighborhood of the diagonal in
$X\times X$.
Let $x\in X$ be arbitrary and choose Heisenberg
coordinates $(z,\theta)$ at $P_0=\pi(x)$ (so that $x$ has coordinates $(0,0)$).
Recalling that $\lambda = a(z)^{-\half } e^{i \theta}$ on $X$,
we have by Lemma
\ref{a2}, 
\begin{equation*} \frac{1}{i}\psi(0,0;z, \theta) = 1 -
\frac{a(0,z)}{ \sqrt{a(z)}} e^{-i\theta}
 = (1-e^{-i\theta}) +
e^{-i\theta}\left [\half |z|^2+O(|z|^3)\right ]\;. \end{equation*}
Thus, $$\Re \left[ \frac{1}{i}\psi(0,0;z, \theta)\right]\ge 
0\quad \mbox{for\ } |\theta|<\frac{\pi}{2},\
|z|<\ep\,,$$ where
$\ep$ is independent of the point $P_0\in M$.
\end{proof}

\subsection{The \szego projector}\label{s-szego}
Having defined $C$, we denote by  $I^*(X \times X, C)$  the space of complex Fourier
integral operators with wave fronts along $C$, as described in \S\ref{s-I*}.    We
define   a   Szeg\"o projector $\Pi$ associated to $\Sigma$ and $C$ to be a
self-adjoint projection $\Pi \in I^0(X \times X, C)$ with principal symbol
$\sigma_{\Pi}$ equal  to the canonical 1/2-density of $C$.   Our purpose now is to
describe  the method of  \cite[\S 4]{BS} \cite{BG} (Appendix) for producing a canonical \szego projector
modulo finite rank operators.  We refer there for further discussion of the method.

 Since by Theorem \ref{oscintth}, $C$ is
parametrized by a function of the form $t\psi(x,y)$, the space
$I^n(X \times X, C)$ consists of Fourier
integral operators of the form  $$ \int_0^{\infty} e^{i t\psi(x,y)} s(x,
y, t) dt,\qquad s \sim \sum_{k = 0}^{\infty} t^{m + n-k} s_k(x,y)\in S^{m+n}(X\times
X\times
\R^+)\;,$$ where
$s_k$ is a smooth function of
$(x,y)$ defined near $x = y$.  
Because $C$ is an equivalence relation,  $I^*(X \times X, C)$  is a $*$-algebra. 
Hence there is an induced algebra product (i.e., a $*$-product) defined on the
amplitudes by 
\begin{equation} I(q_1, \psi) \circ I(q_2, \psi) \sim I(q_1 * q_2, \psi),
\end{equation}
where $\sim$ means roughly that the difference is a smoothing
operator. To be more precise, one  uses the method of stationary phase as  in \cite[\S 4,
(4.14)]{BS}  to determine
the symbol expansion of the  $*$-product $*: I^{n_1} \times I^{n_2} \to I^{n_1+ n_2}$ as a formal infinite series.  We  do not wish to consider its convergence, hence  we will use $\sim$ in the weaker sense that the two sides
agree up to any desired order $I^{-K}$ with a well-defined remainder of lower
order.  The leading order term in the
$*$-product is simply the product of the principal symbols.

The first step in constructing a \szego projector $\Pi$ is to construct an approximate projector
$S$, satisfying 
\begin{equation} \label{S} S^2 \sim S \sim S^*, \;\; S \in I^0(X \times X, C),\;\;\;
s_0(x,x) = \pi^{-m} \det L_X ||d \rho||\;\;\; ( \mbox{mod}\;\; \psi).
\end{equation} The very first term of $s * s \sim s$ gives the equation  (\cite{BS}, (4.10))
\begin{equation} \label{SZ} s_0 = \pi^{-m} (h_{\phi})^{-\frac{1}{2}}\;\;\; \mbox{mod}\;\; \psi, \end{equation}
where $h_{\phi}$ is the (non-vanishing) Hessian determinant in $(\sigma, w)$ of the function $\phi(x, y, w, \sigma) = \psi(x, w) +
\sigma \psi(w, y).$  On the diagonal $x = y$ one has $(h_{\phi})^{\frac{1}{2}} =
(\det L_X) ||d \rho||$.  

We now observe that the amplitude of $S$ is not
unique, since we may add to it any amplitude of the form $(D_t - \psi) a$ where $D_t =
\frac{\partial}{i\partial t}.$ This follows from the fact that $(D_t - \psi) e^{i t \psi}
= 0$ and from the fact that we can integrate by parts in an oscillatory integral.  
We  note that $I(q + (D_t - \psi) a, \psi) = I(q, \psi)$ so the product descends
to the quotient by the ideal $(D_t - \psi) S^M.$ We abbreviate the statement 
$I(q + (D_t - \psi) a, \psi) = I(q, \psi)$ by $q \equiv q + (D_t - \psi) a.$
Since  we may absorb
the terms $\psi (s_0 - \pi^{-m} (h_{\phi})^{-\frac{1}{2}})$ into $t^{m - 1 } s_1$,  the difference
being of the form $(D_t - \psi) a$, we may take $s_0$ to be defined by  (\ref{SZ}). 

We now apply an argument given in \cite{BS} to determine the rest of the amplitude
$s$ from  $s_0$ and by  the equations $s^* = s \sim s * s$.  Indeed, since $S^2
\sim S \sim S^*$, there exists $r \in S^{-1}(X\times X\times\R^+)$ with
\begin{equation} s_0 * s_0 = s_0 + r. \end{equation}
Clearly, $r * s_0 = s_0 * r.$ We now construct $E(r)$ so that 
\begin{equation} s = s_0 + E(r) \implies s * s \sim s, \;\; s^* = s. \end{equation}
Omitting the $*$ since the product is commutative, we get:
\begin{equation}\begin{array}{l} (2 s_0 - 1)  E(r) + E(r)^2 = - r \implies  2
\sqrt{r + \frac{1}{4}} E(r) + E(r)^2 = - r  \\ \\
\implies E(r) \sim \frac{1}{2} -  \sqrt{r + \frac{1}{4}} . \end{array}
\end{equation} The right hand side is a formal power series in $*$ products of
$r$. Given a prescribed symbol order in advance, we may truncate the $*$ product
to obtain a convergent expansion for $E(r)$ plus a remainder which is lower than
the prescribed order. The resulting amplitude
 $s$ will satisfy $s^* = s \sim s * s$ to the prescribed order and hence be   an
amplitude for  an approximate \szego projector $S$.  Thus, for any $K$ we can
construct an approximate 
\szego projector $S_K$ such that $S_K^2 = S_K+ E_K = S_K^* + E_K^*$ where 
$E_K(x,y)$  lies in $\ccal^K(X \times X)$. As discussed in \cite{BS},  one can in
fact produce an approximate projector modulo smoothing operators.

We then modify $S$ by a smoothing operator $E$ to obtain a true projection (see
\cite{BG}, Appendix A.4):
  Since   $S^2 \sim  S$ and $S =S^*$, it follows that    the spectrum of $S$ is
concentrated near $\{0, 1\}$.  Let $U_0,U_1$ be disjoint open sets containing the points of the
spectrum near $0,1$, respectively, and let $F$ be the analytic (locally constant) function on $U_0\cup
U_1$ given by $F(U_0)=\{0\},\ F(U_1)=\{1\}$. Hence
$F(S):= \Pi$ is a true projection.  More concretely, since $S$ is self-adjoint it has
an eigenfunction expansion $S(x, y) = \sum_{n = 1}^{\infty} \la_n f_n(x) f_n(y).$ We may
collect (modulo a finite dimensional ambiguity) the eigenvalues which cluster at $\{0\}$
and those which cluster at $\{1\}$ to obtain $S(x, y) = \sum_{n: \la_n \in U_0} \la_n f_n(x) f_n(y)
+  \sum_{n: \la_n \in U_1} \la_n f_n(x) f_n(y).$  Then $F(S)(x,y) =    \sum_{n: \la_n \in U_1}  f_n(x) f_n(y).$

To summarize the above discussion, a \szego projector $\Pi$ can be written in the form:
\begin{equation}\label{oscint}\begin{array}{c}\Pi (x,y) = S(x,y)+E(x,y)\;,\\[12pt] S(x,y)=
\int_0^{\infty}
e^{i t \psi(x,y)} s(x,y,t ) dt\,, \qquad E(x,y)\in  \ccal^\infty(X \times X)\,,\end{array}
\end{equation} where  $\psi$ is given by Theorem
\ref{oscintth} and 
$s \sim \sum_{k = 0}^{\infty} t^{m -k} s_k(x,y)\in S^m(X\times X\times
\R^+)$ is constructed as above. Although  $\Pi$ is not
unique,  the above construction defines a canonical
choice of 
$\Pi$ modulo smoothing operators.   In the complex case,  the construction
produces the usual \szego projector $\Pi$ onto the kernel of
$\dbar_b$, and (\ref{oscint}) is  the Boutet de
Monvel-Sj\"ostrand oscillatory integral formula for it ({\rm
\cite [Th.~1.5 and
\S 2.c]{BS}}).

\subsection{Almost holomorphic sections}\label{s-almostholo}
In the complex case, a holomorphic section $s$ of $L^N$ lifts to a
$\hat{s}\in \lcal^2_N(X)$ which satisfying $\dbar_b \hat{s} = 0.$  The
operator $\dbar_b$ extends to a complex satisfying $\dbar_b^2 = 0$,
which is a necessary
and sufficient condition for having a maximal family of CR holomorphic
coordinates.
In the non-integrable case  $\dbar_b^2 \not= 0$, and there may be no local
solutions of $\dbar_b f = 0$.  We now quote a result of Boutet-de-Monvel and Guillemin which  replaces $\dbar_b$ with a
pseudodifferential operator  $\bar{D}_0$  so
that
$\bar{D}_0
\Pi =0$. 
Indeed,
Boutet de Monvel \cite{Bou} and Boutet de Monvel - Guillemin \cite{BG}
defined a complex $\bar{D}_j$, which is a good replacement for
$\dbar_b$ in the non-integrable
case.  Their main result is:

\begin{theo}\label{COMPLEX} {\rm (see \cite{BG}, Lemma 14.11 and Theorem A
5.9)} There exists  an $S^1$-invariant  complex of first order
pseudodifferential operators $\bar{D}_j$ over $X$
$$0 \rightarrow \ccal^\infty(\Lambda_b^{0,0})
\ {\buildrel \bar{D}_0 \over\to}\  \ccal^\infty(\Lambda_b^{0,1})
\ {\buildrel \bar{D}_1 \over\to}\  \cdots \ {\buildrel \bar{D}_{m-1}
\over\longrightarrow}\
\ccal^\infty(\Lambda_b^{0,m})\to 0\,,$$ where $\Lambda_b^{0,j}=\La^j
(H^{0,1}X)^*$, such that:

\begin{enumerate}
\item[{\rm i)}]  $\sigma(\bar{D}_j) =
\sigma(\bar{\partial}_b)$ to second order along $\Sigma:=\{(x,r\al_x):x\in X,
r>0\}\subset T^*X$;
\item[{\rm ii)}]  The \szego kernel $\Pi$ is the orthogonal projector onto the
kernel of $\bar{D}_0$; 
\item[{\rm iii)}] $(\bar{D}_0,
\frac{\partial}{\partial \theta})$ is jointly elliptic. 
\end{enumerate}\end{theo}

Let us briefly summarize the
construction of  $\bar{D}_0$ (following \cite[Appendix]{BG}).
We begin with any $S^1$-equivariant symmetric first order pseudodifferential
operator $\bar{D}_0'$ with
principal symbol equal to $\sum_{j = 1}^m \zeta_j \bar{\vartheta}_j.$  Then
$\bar{D}'_0 \Pi$ is of order
$\leq 0$ so one may  find a zeroth order
pseudodifferential
operator $Q_0$ such that $\bar{D}_0' \Pi \sim Q_0 \Pi$ (modulo smoothing
operators).   Then put:
$\bar{D}_0 =
(\bar{D}_0' - Q_0) -
(\bar{D}_0' - Q_0)\Pi$. Clearly, $\bar{D}_0 \Pi = 0$ and
$\sigma({\bar{D}_0})
= \sigma(\bar{D}_0') = \sum_{j = 1}^m \zeta_j \bar{\vartheta}_j.$
The characteristic variety of $\bar{D}_0$ is then equal to $\Sigma$.  Since
$p_{\theta}$ is the
symbol of $\frac{\partial}{\partial \theta}$ and since the system
$\{\sigma_{\bar{D}_0}, p_{\theta}\}$
has no zeros in $T^*X - 0$ it follows that $\{\bar{D}_0,
\frac{\partial}{\partial \theta}\}$ is an elliptic
system.

One can then construct the higher $\bar{D}_j$ recursively so that
$\bar{D}_j
\bar{D}_{j-1} = 0$. We refer to
\cite{BG}, Appendix \S 5, for further details.

\medskip
We refer to the kernel $\hcal^2(X)=\ker \bar D_0
\cap \lcal^2(X)$ as the   Hardy space
 of square-integrable `almost CR functions' on $X$. The
$\lcal^2$ norm is with respect to the inner product
\begin{equation}\label{inner} \langle  F_1, F_2\rangle
=\frac{1}{2\pi}\int_X
F_1\overline{F_2}dV_X\,,\quad F_1,F_2\in\lcal^2(X)\,,\end{equation}
where \begin{equation}\label{dvx}dV_X=\frac{1}{m!}\al\wedge
(d\al)^m=\al\wedge\pi^*dV_M\,.\end{equation}

The $S^1$ action on $X$ commutes
with $\bar{D}_0$; hence $\hcal^2(X) = \bigoplus_{N
=0}^{\infty} \hcal^2_N(X)$ where $\hcal^2_N(X) =
\{ F \in \hcal^2(X): F(r_{\theta}x)
= e^{i
N \theta} F(x) \}$. We denote by $H^0_J(M, L^{ N})$ the space of sections
which
corresponds to $\hcal^2_N(X)$  under
the  map
$s\mapsto
\hat{s}$.  Elements of $H^0_J(M, L^{ N})$ are the {\it almost holomorphic
sections\/} of $L^N$. (Note that products of almost holomorphic sections are
not necessarily almost holomorphic.) We henceforth write $\hat s =s$ and
identify $H^0_J(M, L^{ N})$ with $\hcal^2_N(X)$.  Since $(\bar{D}_0,
\frac{\partial}{\partial \theta})$
is a jointly elliptic system, elements of $H^0_J(M, L^{ N})$ and
$\hcal^2_N(X)$ are smooth.  In many other respects, $H^0_J(M, L^N)$ is
analogous to the space of
holomorphic sections in the complex case.  Subsequent results will bear
this out.

We let $\Pi_N : \lcal^2(X) \rightarrow \hcal^2_N(X)$ denote the
orthogonal
projection.  The  level $N$ Szeg\"o kernel $\Pi_N(x,y)$ is defined by
\begin{equation} \Pi_N F(x) = \int_X \Pi_N(x,y) F(y) dV_X (y)\,,
\quad F\in\lcal^2(X)\,.
\end{equation} It can be given as
\begin{equation}\label{szego}\Pi_N(x,y)=\sum_{j=1}^{d_N}
S_j^N(x)\overline{ S_j^N(y)}\,,\end{equation} where
$S_1^N,\dots,S_{d_N}^N$ form an orthonormal basis of
$\hcal^2_N(X)$.

\begin{rem}
The results stated here use only the $\bar D_0$ term of the complex; its
kernel consists of the spaces of almost holomorphic sections of the powers
$L^N$ of the line bundle $L$, as explained below. The complex $\bar D_j$ was
used by Boutet de Monvel -Guillemin 
\cite[Lemma~14.14]{BG} to show that the dimension of $H^0_J(M, L^N)$ or $
\hcal^2_N(X)$  is given by the Riemann-Roch formula (for
$N$ sufficiently
large). For our results, we need only the leading term of Riemann-Roch,
which we obtain as a consequence of Theorem \ref{tyz}(a). 
(The reader should be warned that the symbol is
described incorrectly in Lemma 14.11 of \cite{BG}.  However, it is correctly
described
in Theorem 5.9 of the Appendix to \cite{BG} and also in \cite{GU}).
\end{rem}

\section{Scaling asymptotics for Szeg\"o kernels}\label{s-neardiag}

In \cite[Theorem~3.1]{BSZ1}, we showed that in the complex case, the
scaled Szeg\"o kernel $\Pi_N$ near the diagonal is asymptotic to the Szeg\"o
kernel $\Pi^\H_1$ of level one for the reduced Heisenberg group, given by
\begin{equation}\label{Heisenberg}\Pi^\H_1(z,\theta;w,\phi) =
\frac{1}{\pi^m} e^{i(\theta-\phi)+i\Im
(z\cdot \bar w)-\half |z-w|^2}= \frac{1}{\pi^m} e^{i(\theta-\phi)+z\cdot
\bar
w-\half(|z|^2+|w|^2)}\,.\end{equation}  The method was to apply the Boutet de
Monvel-Sj\"ostrand oscillatory integral formula
$$\Pi (x,y) \sim
\int_0^{\infty}
e^{i t \psi(x,y)} s(x,y,t ) dt \;\;\; \mbox{mod}\;\; \ccal^\infty(X \times X)
$$ arising from a parametrix construction (recall (\ref{oscint})). 

Our  goal now is to show that the universal asymptotic formula
of
\cite{BSZ1} for the near-diagonal scaled Szeg\"o kernel holds for the
symplectic case (Theorem \ref{neardiag}).  In fact, our description adds some quite useful details
to the formula given in \cite{BSZ1}.

We have shown that \szego kernels can be expressed  in the form  $\Pi(x,y) = S(x,y)
+ E(x,y)$, where
$S$ is the Fourier integral kernel in (\ref{oscint}) and where $E \in \ccal^\infty(X
\times X)$ is the remainder.  We denote the
$N^{\rm th}$ Fourier coefficient of these operators relative to the $S^1$ action by
$\Pi_N = S_N + E_N$. Since $E$ is smooth, we have $E_N(x,y) = O(N^{-\infty})$, where
$O(N^{-\infty})$ denotes a quantity which is  uniformly $O(N^{-k})$ on $X\times X$ for
all positive
$k$. The Fourier coefficients  $S_N$ are given by
\begin{equation}S_N(x,y)  =   \int_0^{2\pi} e^{- i
N \theta}  S( r_{\theta} x,y)
d\theta  =  \int_0^{\infty} \int_0^{2\pi} e^{- i
N \theta}  e^{it  \psi( r_{\theta} x,y)} s(r_{\theta} x,y,t)
d\theta dt \end{equation} where $r_{\theta}$ denotes the $S^1$ action
on $X$. Changing variables $t \mapsto N t$ gives
\begin{equation} S_N(x,y) = N \int_0^{\infty} \int_0^{2\pi}
e^{ i N ( -\theta + t \psi( r_{\theta} x,y))} s(r_{\theta} x,y,
Nt) d\theta dt\,.\end{equation}

We now determine the scaling limit of the Szeg\"o kernel by continuing the argument of
\cite{BSZ1}, and adding  some new details on homogeneities which will be useful in our 
applications.  To describe the scaling limit at a point $x_0\in X$, we 
choose a Heisenberg chart
$\rho:U,0\to X,x_0$ centered at 
$P_0=\pi(x_0)\in M$.  Recall (\S \ref{s-heisenberg}) that choosing $\rho$ is
equivalent to choosing preferred coordinates centered at $P_0$ and a
preferred local frame
$e_L$ at
$P_0$.  We then write the \szego kernel $\Pi_N$ in terms of these
coordinates:
$$\Pi_N^{P_0}(u,\theta;v,\phi)=\Pi_N(\rho(u,\theta),\rho(v,\phi))\,,$$
where the superscript $P_0$ is a reminder that we are using coordinates
centered at $P_0$.  (We remark that the function $\Pi_N^{P_0}$ depends
also on the choice of preferred coordinates and preferred frame, which we
omit from the notation.)  The first term in our asymptotic formula below
says that the $N^{\rm th}$ scaled \szego kernel looks
approximately like the  Szeg\"o kernel of level one
for the
reduced  Heisenberg group (recall (\ref{Heisenberg})):
$$\Pi_N^{P_0}(\frac{u}{\sqrtn},\frac{\theta}{N};
\frac{v}{\sqrtn},\frac{\phi}{N})\approx \Pi^\H_1(u,\theta;v,\phi)=
\frac{1}{\pi^m} e^{i(\theta-\phi)+i\Im
(u\cdot \bar v)-\half |u-v|^2}\,.$$   In the following, we shall denote the 
Taylor series of a  $\ccal^\infty$ function $f$ defined in a
neighborhood of $0\in \R^K$ by
$f \sim f_0 + f_1 + f_2 + \dots$ where $f_j$ is the homogeneous polynomial
part of degree $j$.
We also denote by $R_n ^f \sim f_{n+1} + \cdots$ the remainder term
in the Taylor expansion.

The following is our main result on the scaling asymptotics of the \szego
kernels near the diagonal.  Since the result is of independent interest, we
state our asymptotic formula in a more precise form than is needed for the 
applications in this paper.

\begin{theo} \label{neardiag}
Let $P_0\in M$ and choose a Heisenberg coordinate chart about $P_0$.
Then
$$\begin{array}{l} N^{-m}\Pi_N^{P_0}(\frac{u}{\sqrtn},\frac{\theta}{N};
\frac{v}{\sqrtn},\frac{\phi}{N})\\ \\ \qquad
= \Pi^\H_1(u,\theta;v,\phi)\left[1+ \sum_{r = 1}^{K}
N^{-r/2} b_{r}(P_0,u,v)
+ N^{-(K +1)/2} R_K(P_0,u,v,N)\right]\;,\end{array}$$
where:
\begin{itemize}

\item $ b_{r} = \sum_{\alpha=0}^{2[r/2]}
\sum_{j=0}^{[3r/2]}(\psi_2)^{\alpha}Q_{r,\al,3r-2j}
\,,$ where $Q_{r,\alpha ,d}$ is homogeneous of
degree $d$ and
\begin{equation*} \psi_2(u,v) =
u \cdot\bar{v} - \half(|u|^2 + |v|^2)\,;\end{equation*}
in particular, $b_r$ has only even
homogeneity if $r$ is even, and only odd homogeneity if $r$ is odd;\\

\item $\|R_K(P_0,u,v,N)\|_{\ccal^j(\{|u|\le \rho,\ |v|\le \rho\}}\le
C_{K,j,\rho}$ for $j\ge 0,\,\rho>0$ and $C_{K,j,\rho}$ is independent of the
point
$P_0$ and choice of coordinates.
\end{itemize}
\end{theo}

\begin{proof} Since $\Pi_N=S_N+O(N^{-\infty})$, it suffices to give the asymptotics  of
$S_N$. Hence, we fix $P_0$ and consider the asymptotics of
 \begin{equation}\begin{array}{l}\displaystyle S^{P_0}_N\left(
\frac{u}{\sqrt{N}}, 0;  \frac{v}{\sqrt{N}}, 0\right)
\\[14pt]\displaystyle
\quad\quad=
 N \int_0^{\infty} \int_0^{2\pi}
 e^{ i N \left( -\theta + t\psi(  \frac{u}{\sqrt{N}}, \theta; 
\frac{v}{\sqrt{N}}, 0)\right)} s\big( \frac{u}{\sqrt{N}}, \theta; 
\frac{v}{\sqrt{N}}, 0, Nt\big) d\theta dt \,,\end{array}\label{SN}\end{equation}
where $\psi$ and $s$ are the phase and symbol from (\ref{oscint})
written in terms of the Heisenberg coordinates.

On $X$ we  have $\lambda = a(z)^{-\half } e^{i \phi}$.
 So for $(x,y) = (z, \phi, w, \phi') \in X \times X$, we have by Theorem
 \ref{oscintth},
\begin{equation} \psi(z, \phi, w, \phi') = i \left[1 -
\frac{a(z,w)}{ \sqrt{a(z)} \sqrt{a(w)}} e^{i
(\phi -
\phi')}\right]\;. \end{equation} It follows that
\begin{equation}\begin{array}{l}\displaystyle\psi( 
\frac{u}{\sqrt{N}},
\theta;  \frac{v}{\sqrt{N}}, 0)\\[14pt] \displaystyle \quad\quad=
i\left[1 - \frac{a( \frac{u}{\sqrt{N}},
\frac{{v}}{\sqrt{N}})}{ \sqrt{a( \frac{u}{\sqrt{N}}, 
\frac{ u}{\sqrt{N}})} \sqrt{a( \frac{v}{\sqrt{N}}, 
\frac{ v}{\sqrt{N}})}} e^{i \theta}\right] .
\end{array}\end{equation}

We observe that the asymptotic expansion of a function $f(
\frac{u}{\sqrt{N}}, 
\frac{v}{\sqrt{N}})$ in powers of $N^{-\half}$ is just the
Taylor
expansion of $f$ at  $u = v = 0$.  By Lemma \ref{a2} and the notational
convention
established above, we have
 \begin{equation}\label{taylorh2} a(
\frac{u}{\sqrt{N}}, 
\frac{v}{\sqrt{N}}) = 1 + \frac{1}{N} u\cdot \bar v +
R_3^a ( \frac{u}{\sqrt{N}},
\frac{v}{\sqrt{N}})\,,\qquad R_3^a ( \frac{u}{\sqrt{N}},
\frac{v}{\sqrt{N}})=O(N^{-3/2})\,.\end{equation}

The phase in (\ref{SN})
\begin{equation} 
 \wt\Psi:= t \psi( \frac{u}{\sqrt{N}},\theta;  \frac{
v}{\sqrt{N}}, 0)  -\theta  = it \left[ 1 - \frac{a(
\frac{u}{\sqrt{N}},  \frac{
v}{\sqrt{N}})}{a( \frac{u}{\sqrt{N}},  \frac{
u}{\sqrt{N}})^{\half} a( \frac{v}{\sqrt{N}},  \frac{
v}{\sqrt{N}})^{\half}} e^{i \theta}\right] -  \theta 
\end{equation}
then has the asymptotic $N$-expansion \begin{equation}\label{entirephase}
\wt\Psi=it[ 1
- e^{i \theta}] -  \theta - \frac{it}{N}\psi_2(u,v) e^{i \theta} + t
R_3^\psi(\frac{u}{\sqrtn},\frac{v}{\sqrtn}) e^{i \theta}\,. \end{equation}

We use a smooth partition of unity $\{\rho_1(t),\rho_2(t)\}$
to decompose the integral (\ref{SN}) into one over $0<t<3$
and one over $t>2$: \begin{eqnarray} S_N^{P_0}(
\frac{u}{\sqrt{N}}, 0;  \frac{v}{\sqrt{N}}, 0) &=& I_1+I_2\;,\nonumber\\
I_1&=&N \int_0^3 \int_0^{2\pi}
 e^{ i N \wt\Psi} \rho_1(t)s\big( \frac{u}{\sqrt{N}}, \theta; 
\frac{v}{\sqrt{N}}, 0, Nt\big) d\theta dt\;,\label{I1}\\
I_2&=&N \int_2^{\infty} \int_0^{2\pi}
 e^{ i N \wt\Psi} \rho_2(t)s\big( \frac{u}{\sqrt{N}}, \theta; 
\frac{v}{\sqrt{N}}, 0, Nt\big) d\theta dt\;.\label{I2}\end{eqnarray}

To evaluate $I_1$, we absorb $ (-i\psi_2 + NR_3^\psi)te^{i \theta}$   into
the
amplitude (as in \cite{BSZ1}),
so that  we view  $I_1$
as
an oscillatory integral with phase \begin{equation}\label{phase}\Psi(t,\theta): = it (
1 - e^{i\theta})-
\theta\end{equation} and with amplitude \begin{equation}\label{amplitude}
A(t,\theta;P_0,u,v):=\rho_1(t)e^{ t e^{i \theta} \psi_2(u,v)  +  it e^{i \theta} N
R_3^\psi(
\frac{u}{\sqrt{N}},\frac{
v}{\sqrt{N}})}   s\big( \frac{u}{\sqrt{N}}, \theta; 
\frac{v}{\sqrt{N}}, 0, Nt\big)\,;\end{equation}
i.e., \begin{equation}\label{phase-amplitude} I_1 = N
\int_0^3\int_0^{2\pi} e^{iN\Psi(t,\theta)}A(t,\theta;P_0,u,v)d\theta
dt\end{equation}

We evaluate (\ref{phase-amplitude}) by the method of
stationary
phase as in \cite{BSZ1}.
The phase $\Psi$ is independent of the parameters $(u,v)$ and we have
\begin{equation} \begin{array}{l} \frac{\d}{\d t} \Psi = i (1 - e^{i \theta}
)
\\[8pt]
 \frac{\d}{\d \theta}  \Psi = t e^{i \theta} - 1 \end{array} \end{equation}
so
the critical set of the phase is the point $ \{t=1, \theta = 0\}$.
The  Hessian $\Psi''$ on the critical set equals
$$ \left( \begin{array}{ll} 0 & 1  \\
1 & i  \end{array} \right)$$
so the phase is non-degenerate and the Hessian operator $L_{\Psi}$ is given
by
$$L_{\Psi}= \langle \Psi''(1, 0)^{-1} D,D\rangle =  2
\frac{\partial^2}{\partial t
\partial \theta} - i
 \frac{\partial^2}{
\partial t^2}\,.$$  
By the stationary phase
method for complex oscillatory integrals (\cite{H}, Theorem
7.7.5), we have
\begin{equation}\label{MSP} I_1 = \ga\sum_{j = 0}^{J} N^{ -j}   L_j [
A(t,\theta;P_0,u,v)]|_{t = 1, \theta = 0} + \wh R_{J}(P_0, u, v, N),\end{equation} where
$$\ga= N \frac{1}{\sqrt{\rm{det} (N \Psi''(1,0)/2
\pi i)}}
 =\sqrt{-2\pi i}$$ and $L_j$ is the differential operator of order $2j$ in
$(t,\theta)$ defined by
\begin{equation} L_j \phi(t,\theta)  =
\sum_{\nu - \mu = j}\sum_{2
\nu \geq 3 \mu}
\frac{1}{2^{\nu} i^j \mu! \nu!}
 L_{\Psi}^{\nu} [ \phi(t,\theta)  (R_3^\Psi)^{\mu}(t, \theta)]
\end{equation}
with $R_3^\Psi(t,\theta)$  the third order remainder in the Taylor
expansion
of
$\Psi$ at $(t,\theta) = (1,0)$. 

Also, the remainder is estimated by
\begin{equation} \label{remainder} |\wh R_{J}(P_0, u, v, N)| \leq C
N^{- J}
\sum_{|\alpha| \leq 2J+2} \sup_{t, \theta} |D^{\alpha}_{t,
\theta} A(t,\theta;P_0,u,v) |\le C_JN^{m-J}\;,
\end{equation} where the last inequality follows by observing from (\ref{amplitude}) that
\begin{equation}|A(t,\theta;P_0,u,v)|\le C' |s\big( \frac{u}{\sqrt{N}}, \theta; 
\frac{v}{\sqrt{N}}, 0, Nt\big)|\le C''N^m\qquad (0\le t\le
3)\;,\label{A0}\end{equation} and similarly for its derivatives (using the fact
that $s$ is a symbol of order $m$ and hence $D^{\alpha}_{t,
\theta} s(r_\theta x,y,Nt)=O(N^m)$ uniformly for $t\le 3$).

To evaluate (\ref{MSP}), we first expand $\exp\left[ it e^{i \theta} N
R_3^\psi(\frac{u}{\sqrt{N}},\frac{
v}{\sqrt{N}})\right]$ in powers of $N^{-\half}$, keeping track of the
homogeneity in $(u,v)$ of the
coefficients.  We simplify the notation by writing $g= t e^{i
\theta}$.  By definition,
$$R_3^\psi( \frac{u}{\sqrt{N}}, 
\frac{v}{\sqrt{N}}) \sim N^{-3/2} \psi_3(u,v) +  N^{-2} \psi_4(u,v) + \cdots
+ N^{-d/2} \psi_d(u,v) + \cdots.$$
 We then  have
\begin{equation}\label{EXPAND} e^{  iN g R_3^\psi(
\frac{u}{\sqrt{N}}, 
\frac{
v}{\sqrt{N}})}\sim \sum_{r = 0}^{\infty} N^{-r/2}
c_{r}(u,v;t,\theta)\,.\end{equation} 

We further expand
\begin{equation}\label{EXPAND2}
c_{r} = \sum_{\la = 1 }^{r} c_{r,r+2\la}(u,v; t,\theta)\,,\ r\ge 1\,,\quad
c_0 = c_{00}=1\,,\end{equation}
with $c_{rd}$ homogeneous of degree $d$ in
$u,v$.  Note that $c_r$ is a
polynomial of degree $r$ in
$g$. 
(The explicit formula for $c_{rd}$ is:
$$\textstyle c_{rd} = \sum \left\{ \frac{1}{n!}(ig)^n
\Pi_{j=1}^n
\psi_{a_j}(u,v):
n\ge 1,  a_j \geq 3, \sum_{j = 1}^n a_j=d,
\sum_{j=1}^n (a_j - 2)=r
\right\}\,,\ r\ge 1\,.$$  The range of $d$ is determined by the fact
that
$d =\sum_{j = 1}^n a_j = r + 2n$ with $0 \leq n \leq r$.)

We similarly expand the symbol:  $$s(\frac{u}{\sqrt{N}},
\theta;  \frac{ v}{\sqrt{N}},0,Nt)=\sum_{k =
0}^{\infty} N^{m-k} t^{m-k} s_k(
\frac{u}{\sqrt{N}},  \frac{
v}{\sqrt{N}}, \theta)=
\sum_{k, \ell = 0}^{\infty} N^{m-k- \ell/2} t^{m-k} s_{k \ell}( u,v,\theta) $$
where $s_{k \ell}$ is the homogeneous term of $s_k$ of degree $\ell$ in $(u,v)$.
Hence, we have
\begin{eqnarray}A &\sim &  \rho_1(t)e^{ g\psi_2(u,v)+iNg R_3^\psi(
\frac{u}{\sqrt{N}}, 
\frac{
v}{\sqrt{N}}, \theta)} \sum_{k = 0}^{\infty} N^{m-k} t^{m-k} s_k(
\frac{u}{\sqrt{N}},  \frac{v}{\sqrt{N}}, \theta)\nonumber\\& \sim&  \rho_1(t) e^
{g\psi_2(u,v)}
 N^m \sum_{n = 0}^{\infty} N^{-n/2} f_{n}(u, v;t, \theta,  P_0)
\;, \label{A}\end{eqnarray}
where  the remainder in (\ref{A}) after summing $K$ terms is
$O(N^{m-\frac{K+1}{2}})$. (Note that
$f_n$ is a polynomial in $t$ of degree
$n+m$.) 
We further have
 $$f_{n}
= \sum_{r + \ell + 2k  = n} c_{r} s_{k \ell}
= \sum_{k=0}^{[n/2]}  t^{m - k}\left( s_{k,n-2k}
+\sum_{r=1}^{n-2k}\sum_{\la=1}^{r} c_{r,r+2\la}  s_{k,n-2k-r}\right)
=\sum_{j=0}^{[3n/2]}f_{n,3n-2j}\;,$$
where $f_{n,d}$ is homogeneous
of degree $d$ in $(u,v)$.

Hence, \begin{equation}\label{LjA}
L_j[A]|_{t = 1, \theta = 0}=\sum_{n=0}^K N^{m-n/2}L_j [ e^{g \psi_2}
f_{n}]|_{t = 1,
\theta = 0}+O(N^{m-\frac{K+1}{2}})\;.\end{equation}
Since $L_{\Psi}$ is a second order operator in  $(t,
\theta)$, we see that
\begin{equation} \label{SERIES} 
L_j [ e^{g \psi_2} f_{n}]|_{t = 1, \theta = 0}= e^{\psi_2}
\sum_{\alpha \leq 2j} (\psi_2)^{\alpha} F_{n j \alpha}\,,
\end{equation} where the $F_{nj\al}$ are polynomials in $u,v,\bar u,\bar v$ of degree $\le 3n$.
Therefore, by (\ref{MSP}), (\ref{LjA}) and (\ref{SERIES}), we have
\begin{eqnarray}\label{MSP2}N^{-m}I_1 &\sim
&\sqrt{-2\pi i}\,e^{\psi_2}\sum_{n=0}^\infty \sum_{j=0}^\infty \sum_{\al=0}^{2j}(\psi_2)^\al
N^{-\frac{n}{2}-j} F_{nj\al}\nonumber\\ &\sim
&\sqrt{-2\pi i}\,e^{\psi_2} \sum_{r=0}^\infty
\sum_{j=0}^{[r/2]}\sum_{\al=0}^{2j}(\psi_2)^\al
N^{-r/2}F_{r-2j,j,\al}\nonumber\\&\sim
&e^{\psi_2} \sum_{r=0}^\infty \sum_{\al=0}^{2[r/2]}(\psi_2)^\al
N^{-r/2}Q_{r\al}\,.\end{eqnarray}

As with $f_{n}$ we have the homogeneous expansion:
\begin{equation} \label{Q} Q_{r\al}=\sum_{j=0}^{[3r/2]}Q_{r,\al,3r-2j}\,.
\end{equation}
Here, $Q_{r,\al,d}$ is homogeneous of degree $d$ in $(u,v)$. 
Thus we have the desired Taylor series for $I_1$. 

To show that $I_2= O(N^{-\infty})$, we observe from (\ref{entirephase}) that
$$\left|\frac{\d}{\d\theta}\wt\Psi\right|\ge \left[1-O(N\inv)  
\right]t-1\ge \frac{1}{3}t \quad\mbox{for}\ \ t\ge 2,\ N\gg 1\,.$$
Hence by \cite{H}, Theorem 7.7.1, for all integers $k\ge 1$ we have
\begin{eqnarray*} \int_0^{2\pi}
 e^{ i N \wt\Psi} \rho_2(t)s\big( \frac{u}{\sqrt{N}}, \theta; 
\frac{v}{\sqrt{N}}, 0, Nt\big) d\theta \hspace{-1in}& \\ &\le &
 C' N^{-k} t^{-k} \sum_{\al\le k}\sup \left|D^\al_\theta\left(\rho_2(t)s\big(
\frac{u}{\sqrt{N}}, \theta; 
\frac{v}{\sqrt{N}}, 0, Nt\big)\right)\right|\\ & \le & C'_k  N^{m-k}
t^{m-k}\;.\end{eqnarray*}  
(To see that $C'$ is independent of $t$, we write $iN\wt\Psi= i(Nt)\Phi$, with
$\Phi=\wt\Psi/t$ and apply \cite[Theorem 7.7.1]{H} to the phase $\Phi$.) Integrating
over
$t$, we conclude that
$I_2=O(N^{m+1-k})$;
since $k\ge 1$ is arbitrary, $I_2=O(N^{-\infty})$.  The estimate for the remainder  now follows
from  (\ref{remainder}) and (\ref{LjA}). 
\end{proof}

\section{Kodaira embedding and Tian almost isometry
theorem}\label{s-kodaira}

\begin{defin} By the Kodaira maps we mean the maps
$\Phi_N : M \to PH^0(M,L^N)'$ defined by $\Phi_N(z) = \{s^N: s^N(z) =
0\}$. Equivalently, we can choose an orthonormal basis
$S^N_1,\dots,S^N_{d_N}$
of $H^0(M,L^N)$ and write
\begin{equation}\label{Kmap} \Phi_N : M \to\CP^{d_N-1}\,,\qquad
\Phi_N(z)=\big(S^N_1(z):\dots:S^N_{d_N}(z)\big)\,.\end{equation}
We also define the lifts of the Kodaira maps:
\begin{equation}\label{lift}\wt{\Phi}_N : X \to
\C^{d_N}\,,\qquad
\wt\Phi_N(x)=(S^N_1(x),\dots,S^N_{d_N}(x))\,.\end{equation}
\end{defin}

\medskip Note that
\begin{equation}\label{PiPhi} \Pi_N(x,y)=\wt\Phi_N(x) \cdot
\overline{\wt\Phi_N(y)}\,;\end{equation} in particular,
\begin{equation}\label{PiPhi2}
\Pi_N(x,x)=\|\wt\Phi_N(x)\|^2\,.\end{equation}

We now prove Theorem \ref{tyz}, which generalizes to the symplectic category 
the asymptotic expansion
theorem of \cite{Ze} (also proved independently by \cite{Cat} using the
Bergman kernel in place of the \szego kernel) and Tian's approximate
isometry theorem
\cite{Ti}:

\smallskip\noindent (a) Using the expansion of Theorem
\ref{neardiag} with $u=v=0$ and noting that $b_r(z,0,0)=0$ for $r$ odd,
we obtain the above expansion of $ \Pi_N(z,0;z,0)$ with
$a_r(z)=b_{2r}(z,0,0)$. (The expansion also follows by precisely
the same proof as in
\cite{Ze}.)

\smallskip\noindent (b)  In the holomorphic case, (b) followed by
differentiating (a), using
that $\Phi_N^* ( {\partial} \bar\partial \log |\xi|^2) =  {\partial}
\bar\partial
\log |\Phi_N|^2$.  In the almost complex case, $\Phi_N^*$ does
not commute with
the complex derivatives, so we need to modify the proof.  To do so, we
use the following notation: the exterior derivative on a product
manifold $Y_1\times Y_2$ can be decomposed as $d=d^1+d^2$, where
$d^1$ and
$d^2$ denote exterior differentiation on the first and second factors,
respectively. (This is formally analogous to the decomposition $d=\d +\dbar$;
e.g., $d^1d^1=d^2d^2=d^1d^2+d^2d^1=0$.) 

Recall that the Fubini-Study form $\omega_{FS}$  on
$\CP^{m-1}$ is induced by the
2-form  $\wt\omega_m= \frac{i}{2} \ddbar\log |\xi|^2$ on $\C^m\sm\{0\}$.
We consider the 2-form $\Om$ on
$(\C^m\sm\{0\})\times(\C^m\sm\{0\})$ given by
$$\Om=\frac{i}{2}\ddbar \log \zeta\cdot\bar\eta=\frac{i}{2}d^1 d^2 \log
\zeta\cdot\bar\eta\,.$$  Note that $\Om$ is smooth on a neighborhood of the
diagonal $\{\zeta=\eta\}$, and
$$\Om|_{\zeta=\eta}=\wt\om_m$$ (where the restriction to $\{\zeta=\eta\}$
means the pull-back under the map $\zeta\mapsto (\zeta,\zeta)$).

It suffices to show that
$$\frac{1}{N} \tilde{\Phi}_N^* \omega_{d_N} \to \pi^* \omega, \;\;\;\; \pi:
X \to M.$$  To do this, we consider the maps
$$\Psi_N=\wt\Phi_N\times\wt\Phi_N:X\times X \to \C^{d_N}\times
\C^{d_N}\,,\quad
\Psi_N(x,y)=  (\wt\Phi_N(x),\wt\Phi_N(y))\,.$$
It is
elementary to check that
$\Psi_N^*$ commutes with $d^1$ and $d^2$.
By (\ref{PiPhi}),
we have
$$\Psi_N^*(\log \zeta\cdot\bar\eta) = (\log \zeta\cdot\bar\eta)\circ\Psi_N
=\log\Pi_N\,.$$

Therefore, \begin{equation}\label{PsiN}\frac{1}{N} \tilde{\Psi}_N^*
\Omega_{d_N}=
\frac{i}{2N}{\Psi}_N^*d^1 d^2 \log \zeta\cdot\bar\eta
=\frac{i}{2N} d^1d^2 {\Psi}_N^* \log \zeta\cdot\bar\eta
=\frac{i}{2N} d^1d^2 \log\Pi_N\,.\end{equation}

Restricting (\ref{PsiN}) to the diagonal, we then have
$$\frac{1}{N} \tilde{\Phi}_N^* \omega_{d_N} =
\frac{i}{2N} (d^1d^2 \log\Pi_N)|_{x=y}= \diag^*(d^1d^2 \log\Pi_N)\,,$$
where $\diag:X\to X\times X$ is the diagonal map $\diag(x)=(x,x)$.

Using Heisenberg coordinates as in Theorem \ref{neardiag}, we have by the
near-diagonal scaling asymptotics
\begin{eqnarray*} \left.\frac{1}{N} \tilde{\Phi}_N^*
\omega_{d_N}\right|_{P_0}&=& \left.\frac{i}{2N} \diag^*d^1d^2
\log\Pi_N^{P_0}(\frac{u}{\sqrtn},\frac{\theta}{N};
\frac{v}{\sqrtn},\frac{\phi}{N})\right|_0\\
&=&\left.\frac{i}{2N}\diag^* d^1d^2
\log\Pi_1^\H(u,\theta;v,\phi)\right|_0 +
O(N^{-\half})\,.\end{eqnarray*} Finally,
\begin{eqnarray} \left.
\frac{i}{2N}\diag^* d^1d^2
\log\Pi_1^\H(u,\theta;v,\phi)\right|_0 &=&\frac{i}{2N}
\diag^*d^1d^2\big[ i(\theta-\phi) + u\cdot\bar v -
\half(|u|^2+|v|^2)\big]\nonumber \\
&=& \frac{i}{2N} \sum_{q=1}^{m}du_q\wedge d\bar u_q \ = \
\frac{i}{2} \sum_{q=1}^{m}dz_q\wedge d\bar z_q
\ = \ \om|_{P_0}.\label{dxdy}\end{eqnarray}
\qed

\begin{rem}  A more explicit way to show (b) is to expand the Fubini-Study
form:
$$\wt\omega_m=
\frac{i}{2}|\xi|^{-4}\left[ |\xi|^2 \sum_{j =1}^m d\xi_j \wedge
d\bar{\xi}_j -\sum_{j, k = 1}^m \bar{\xi}_j \xi_k d\xi_j \wedge
d\bar{\xi}_k \right]\,.$$ Then
$$\frac{1}{N} \tilde{\Phi}_N^* \omega_{d_N} =
\frac{i}{2}\Pi_N(x,x)^{-2}\{ (\Pi_N(x,x) d^1 d^2 \Pi_N(x,y) -
d^1 \Pi_N(x,y) \wedge d^2 \Pi_N(x,y)\}|_{x = y}\,,$$
and (b) follows from a short computation using Theorem
\ref{neardiag} as above.
\end{rem}

\medskip
It follows from Theorem \ref{tyz}(b) that  $\Phi_N$ is an immersion for
$N\gg
0$.  Using in part an idea
of Bouche \cite{Bouche}, we now give a simple proof that the Kodaira embedding
theorem holds
for symplectic manifolds (Theorem \ref{kodaira}):

\smallskip\noindent Let $\{P_N, Q_N\}$ be any sequence of distinct points
such that
$\Phi_N(P_N) = \Phi_N(Q_N)$.  By passing to a subsequence we may assume that
one of the following two cases holds:

\begin{enumerate}

\item [(i)]  The distance  $r_N:= \mbox{\rm dist}(P_N, Q_N)$ between
$P_N, Q_N$ satisfies
$r_N \sqrt{N} \to \infty;$

\item [(ii)]  There exists a constant $C$ independent of $N$ such that
$r_N\leq C/ \sqrt{N}.$
\end{enumerate}

\smallskip
To prove that case (i) cannot occur, we let $\Pi_N^{P_N}(x) = \Pi_N(x, P_N)$
denote the `peak section' at $P_N$.
By Theorem \ref{neardiag}, we have
$$N^{-m}\int_{B(P_N, r_N)} |\Pi_N^{P_N} |^2 dV \geq 1 - o(1)\;.$$
The same inequality holds for $Q_N$.  If $\Phi_N(P_N) = \Phi_N(Q_N)$ then
the
total $\lcal^2$-norm of $\Pi_N(x, \cdot)$
would have to be $\sim 2 N^{m/2}$, contradicting the asymptotic  $\sim N^{m/2}$
from Theorem \ref{tyz}(a).

To prove that case (ii) cannot occur, we assume on the contrary that
$\Phi_N(P_N)=
\Phi_N(Q_N)$, where $P_N=\rho_N(0)$ and $ Q_N
=\rho_N(\frac{v_N}{\sqrtn})$, $0\ne |v_N|\le C$, using a Heisenberg
coordinate chart $\rho_N$ about
$P_N$. We consider the function
\begin{equation}\label{fN}
f_N(t)=\frac{|\Pi_N^{P_N}(0,\frac{tv_N}{\sqrtn})|^2}{\Pi_N^{P_N}(0,0)
\Pi_N^{P_N}(\frac{tv_N}{\sqrtn},\frac{tv_N}{\sqrtn})}\,.\end{equation}
Then $f_N(0)=1$ and we see from  (\ref{PiPhi}) and the
Cauchy-Schwartz inequality that $f_N\le 1$.
Furthermore, since  we are supposing that $\Phi_N(P_N)=
\Phi_N(Q_N)$, we also have $f_N(1)=1$. Thus for some value of $t_N$ in the
open interval $(0,1)$, we have $f_N''(t_N)=0$.  By Theorem
\ref{neardiag},
\begin{equation}\label{fNexpand}f_N(t)=e^{-|v_N|^2t^2}\left[1+N^{-1/2}\wt
R_{N}(tv_N)\right]\,,\end{equation} where
$$\wt R_{N}(v)=R_1(P_N;0,v,N)+R_1(P_N;v,0,N) - R_1(P_N;v,v,N)
-R_1(P_N;0,0,N) +O(N^{-1/2}) \,.$$
The estimate for $R_1$ yields:
\begin{equation}\label{fNerror} \|\wt R_{N}\|_{\ccal^2\{|v|\le C\}}=
O(1)\,\end{equation}
Since $f_N(1)=1$, it follows from (\ref{fNexpand})--(\ref{fNerror}) that
$|v_N|^2=O(N^{-1/2})$.  (A more careful analysis shows that we can replace
$N^{-1/2}$ with $N^{-1}$ in (\ref{fNexpand}) and thus $|v_N|=O(N^{-1/2})$.)

Write $e^x=1+x+x^2\phi(x)$.  We then have
$$f_N(t)=1-|v_N|^2t^2+|v_N|^4 t^4 \phi(|v_N|^2t^2)+ N^{-1/2}\wt R_N(tv_N)
\left[1-|v_N|^2t^2+|v_N|^4 t^4 \phi(|v_N|^2t^2)\right]\,.$$
Thus by (\ref{fNerror}), $$ f_N''(t)= -2|v_N|^2 +O(|v_N|^4)
+O(N^{-1/2}|v_N|^2)\,,\quad |t|\le 1\,.$$  Since $|v_N|=o(1)$, it follows
that $$0=f_N''(t_N)=(-2+o(1))|v_N|^2\,,$$ which contradicts the assumption
that $v_N\ne 0$.
\qed

\section{Transversal sections}\label{s-transversal}

\subsection {Quantitative transversality}
\label{s-quanttransverse} 

To illustrate  the connection between our almost holomorphic sections and the
asymptotically holomorphic sections of Donaldson and of Auroux, we now explain how
to adapt Donaldson's proof of the existence of quantitatively transverse
asymptotically holomorphic sections to our setting of almost holomorphic sections.
We shall show that the normalized coherent states in our spaces $H^0_J(M, L^N)$ of
almost holomorphic sections satisfy the estimates of
Proposition 34 in \cite{DON.1}, suitably modified.  Using these modified estimates
together with the proof of
Proposition 15 and of Lemma 14 in  \cite{DON.1}, one obtains the following
result on the existence of almost holomorphic transversal sections. 

\begin{prop} There is a constant $C<+\infty$ such that 
for all sufficiently large $N$, there exists
a section
$s_N \in H^0_J(M, L^N)$ with  $\|\dbar s_N\|_\infty \le C$ and which satisfies
$|\partial s_N| >
{C}\inv{\sqrt{N}}$ on $Z_{s_N}.$\label{transversal} \end{prop}

Here $Z_{s_N}$ denotes the zero set of the section ${s_N}$.
The existence of asymptotically holomorphic sections satisfying these
transversality estimates  was given by Donaldson \cite[Theorem~5]{DON.1}.
(Donaldson expresses transversality in terms of the derivatives $\d_{\rm
scaled}=\frac{1}{\sqrtn}\d$,
$\dbar_{\rm scaled}=\frac{1}{\sqrtn}\dbar$ with respect to the scaled coordinates.)
Recall that transversality implies that 
$Z_{s_N}$ is a symplectic submanifold of $M$ (whose homology class is the Poincar\'e
dual of $N[\om]$).

Before we begin the proof of Proposition \ref{transversal}, we recall that the
starting point of Donaldson's construction is the existence  of  concentrated
asymptotically holomorphic sections
$\sigma_p^N$ supported on a ball of radius $O(N^{-1/3})$ about each point $p \in M$
satisfying certain $\ccal^2$ estimates \cite[Proposition 11]{DON.1}.  In the
integrable case, Donaldson constructs {\it holomorphic\/} sections $\sigma_p^N$
satisfying the estimates (\cite{DON.1}, Proposition 34)
$$\begin{array}{ll} e^{-b  d_N(p, q)^2} \leq |\sigma_p^N(q)|  \leq e^{- a 
d_N(q,p)^2}, \  &  d_N(q,p) \leq \epsilon N^{1/6} \\ & \\ |\sigma_p^N(q)|  \leq
e^{- a  N^{1/3}},&  d_N(q,p) \geq \epsilon N^{1/6}
\end{array}$$
where $d_N(q,p) = \sqrt{N} d(q,p)$ denotes the scaled distance on $M$.
Given such sections, one forms the 
$\frac{1}{\sqrt{N}}$-`lattice' of points
$\Lambda_N = \{p_i\}$ and  the associated  complex vector space of sections
\begin{equation} \label{LC} s_N = \sum_i w_i
\sigma_{p_i}^N, \;\;\;\; w_i \in \C\,. \end{equation}
Donaldson  shows (\cite{DON.1}, Proposition~32) that one can choose
the coefficients $\{w_i\}$ with $|w_i|< 1$ so that the holomorphic section $s_N$
satisfies the quantitative transversality estimate
\begin{equation} \label{QT} |{\partial s_N}(z)| > {C}\inv{\sqrt{N}},\;\;\; \forall
z\in Z_{s_N}.
\end{equation}

Our principal claim is that the \szego kernel $\Pi_N(\bullet,p)$ itself, when 
normalized to have modulus one at $z = p$, satisfies estimates sufficient to prove
the existence of transversal almost-holomorphic sections of the form (\ref{LC}). 
We define 
\begin{equation}\label{coherent} \sigma^N_p(x) = \frac{\Pi_N(x,p)}{\Pi_N(p,p)}\;,
\end{equation}
so that $\sigma^N_p(p)=1$.
(Recall that
$\Pi_N(p,p) = \pi^{-m}N^m +O(N^{m-1})$.)

\begin{lem}  \label{DECAY} For all $\ep,D\in\R^+$, we have for $N$ sufficiently
large:
$$\begin{array}{rll} {\rm(i)}  &|\sigma_p^N(q)|
\ge [1-\frac{C}{\sqrtn} d_N(q,p)]e^{- \frac{1+\ep}{2} d_N(q,p)^2},
&   d_N(q,p)
 \leq  D\;; \\ & & \\
 {\rm(ii)}  &|\sigma_p^N(q)|
\leqsim [1+\frac{C}{\sqrtn} d_N(q,p)]e^{- \frac{1-\ep}{2} d_N(q,p)^2},
&  d_N(q,p)
 \leq  N^{1/6}\;; \\ & & \\
{\rm(iii)}  & N^{-k/2}|\nabla^k\sigma_p^N(q)| \leqsim C_k e^{-\frac{1-\ep}{2}
d_N(q,p)^2}
\;,
&  k\ge 0; \\ & & \\
{\rm(iv)} & N^{-k/2}|\nabla^k\dbar\sigma_p^N(q)| \leqsim C_k e^{-\frac{1-\ep}{2}
d_N(q,p)^2}
\;,
&  k\ge 0; \;
\end{array}$$
where $\al\leqsim\be$ means that  $\al-\be\le O(N^{-\infty})$, uniformly for 
$(q,p)
\in M \times M$. \end{lem}  
\noindent Since the length  $|\sigma_N(p,\theta)|$
is independent of $\theta$, we have dropped $\theta$ from the notation and
regard 
$|\sigma_N|$ as a function on
$M$.  Note that (iii) implies that
\begin{equation}|\nabla^k\sigma_p^N(q)| = O(N^{-\infty})\,,\quad \rm{for}\  \ d_N(q,p)
\ge  N^{1/6},\quad k=0,1,2,\dots\;.\label{outsidebigball}\end{equation}

\begin{proof} We let $q=p+\frac{u}{\sqrtn}$ (in preferred coordinates), so that
$u\approx d_N(q,p)$,  and we write
$\sigma_p^N(q)=\sigma_p^N(\frac{u}{\sqrtn},0)$.  The lower bound (i) is an
immediate consequence of Theorem \ref{neardiag}. To verify (ii)--(iv)
for
$u\le N^{1/6}$ (i.e., on the balls of unscaled radii
$ N^{-1/3}$), we observe that the asymptotic expansion of Theorem
\ref{neardiag} has the following extension to the scaled
$N^{1/6}$ balls:
$$\begin{array}{l} N^{-m}\Pi_N^{P_0}(\frac{u}{\sqrtn},\frac{\theta}{N};
0,0)\\ \\ \qquad
= \Pi^\H_1(u,\theta;0,0)\left[1+ \sum_{r = 1}^{K}
N^{-r/2} b_{r}(P_0,u,0)\right]
+ N^{-(K +1)/2} R_K(P_0,u,N)\;,\end{array}$$ where
\begin{equation}\label{better}\|R_K(P_0,u,N)\|_{\ccal^j(\{|u|\le
N^{1/6}\}}\leqsim C_{K,j}e^{-\frac{1-\ep}{2} |u|^2} \quad \mbox{for}\ j\ge
0\,.\end{equation}
To verify (\ref{better}), we modify the proof of Theorem \ref{neardiag} as follows:
We first note that
\begin{equation} NR_3^\psi( \frac{u}{\sqrt{N}},0 )
 \le C\;,\quad \mbox{for}\ \ |u|\le N^{1/6}\,,\label{R3}\end{equation}
and hence we obtain the large-ball version of (\ref{A0}):
\begin{equation}\label{A1}\begin{array}{r}|A(t,\theta;P_0,u,0)|\le C'e^{-\Re
g|u|^2/2} |s\big(
\frac{u}{\sqrt{N}},
\theta; 
0, 0, Nt\big)|\le C''e^{-\frac{1-\ep}{2}|u|^2}N^m\;,\\[8pt]  1-\de\le t\le
3\;,\ -\de\le\theta\le\de\;,\
|u|\le 2N^{1/6}\;.\end{array}\end{equation}  We use a smooth partition of unity to
decompose the integral (\ref{phase-amplitude}):
\begin{equation}\label{phase-amplitude1} I_1 = I'_1+I''_1,\qquad I'_1= N
\int_{1-\de}^3\int_{-\de}^{\de} e^{iN\Psi(t,\theta)}A(t,\theta;P_0,u,0)d\theta
dt\;.\end{equation} The integral $I''_1$ is over a compact region outside
the critical point
$(1,0)$ of the phase $\Psi$.  In fact $\left|\frac{\d\Psi}{\d\theta}\right|>\al>0$
in this region, and the method of proof that $I_2=O(N^{-\infty})$ shows that we also
have
$I''_1=O(N^{-\infty})$, uniformly for $|u|\le N^{1/6}$ by (\ref{A1}). 

Furthermore by (\ref{R3}), 
the asymptotic expansion (\ref{EXPAND}) of the exponential $e^{N g
R_3(\frac{u}{\sqrt{N}},0 )}$ holds for $|u|\le N^{1/6}$, and hence we have the
large-ball version of (\ref{A}):
\begin{eqnarray} A(t,\theta;P_0,u,0) &=&  \rho_1(t) e^
{-g|u|^2/2}
 N^m \left[\sum_{n = 0}^{K} N^{-n/2} f_{n}(u, 0;t, \theta, 
P_0)+R_K(u,t,\theta)\right]\,,\nonumber\\ &&\qquad\qquad\qquad|R_K|\ \le \
C_KO\big(N^{-\frac{K+1}{2}}\big)\,.
\label{Alarge}\end{eqnarray}
The remainder estimate (\ref{better}) now follows from (\ref{A1}) and (\ref{Alarge})
as in the proof of Theorem \ref{neardiag}. 

To verify (ii) for $|u|\le N^{1/6}$, we apply (\ref{better}) (with $K=0$) to
conclude that
\begin{equation} \label{neardiag2}\textstyle |\sigma^N_p
(p+\frac{u}{\sqrtn})| =
e^{-\half|u|^2}+\frac{1}{\sqrtn}e^{-\frac{1-\ep}{2}|u|^2}R(p,u,N)\,,\quad  |D^\al_u
R(p,u,N)|\le C_{|\al|}\;.\end{equation}
Since
$|\sigma^N_p(p)|=1$,  it follows that $R(p,0,N)=0$ and hence
\begin{equation}\label{R}R(p,u,N)\le C|u|\,,\quad\mbox{for}\ \ 
|u|<N^{1/6}\;.\end{equation} 
Since
$|u|
\approx d_N(p,p+\frac{u}{\sqrtn})$ for $|u|<\ep \sqrtn$, the inequality (ii)
follows  from (\ref{neardiag2})--(\ref{R}). 

To show that (iii)--(iv) holds for $  d_N(q,p) \le  N^{1/6}$, we lift a local
frame
$\{\bar Z_j^M\}$ of the form (\ref{localtangentframe}) to obtain the local frame
$\{\bar Z_1,\dots,\bar Z_m\}$ for
$H^{0,1}X$ given by
\begin{equation}\label{barz}\bar Z_j(p+\frac{u}{\sqrtn}) = \frac{\d^h}{\d\bar z_j} +
\sum_{r=1}^m B_{jk}(\frac{u}{\sqrtn})\frac{\d^h}{\d z_k}\,,\qquad
B_{jk}(\frac{u}{\sqrtn})=O(\frac{1}{\sqrtn}|u|)\le O(N^{-1/3})\,.\end{equation} 
Recalling (\ref{dhdzj}) and (\ref{dhoriz}), we further have
\begin{equation}
\frac{\d^h}{\d\bar z_j} = 
\frac{\d}{\d\bar z_j}
+\left[ -\frac{i}{2}\frac{u_j}{\sqrtn}-R_1^{\bar A_j}
(\frac{u}{\sqrtn})\right]\frac{\d}{\d \theta}
= \frac{\d}{\d\bar z_j}
+i\sqrtn \,u_j-
iNR_1^{\bar A_j}(\frac{u}{\sqrtn})\quad \mbox{at}\ \ 
p+\frac{u}{\sqrtn}\,,\label{dh}\end{equation}  where
$R_1^{\bar A_j}(\frac{u}{\sqrtn})=O(\frac{1}{N}|u|^2)$. The
inequalities (iii)--(iv)  for $|u| \le  N^{1/6}$ follow as before, using
(\ref{barz})--(\ref{dh}).

It remains to show (\ref{outsidebigball}).  For this it
suffices to show that
$\nabla^k\Pi_N^p(\frac{u}{\sqrtn},0;0,0) = O(N^{-\infty})$ for $|u|\ge N^{1/6}$. Since
$\Pi_N = S_N +E_N$, where
$E_N(x,y) = O(N^{-\infty})$, we can neglect the remainder $E_N$. Furthermore, since
$S_N$ is rapidly decaying away from the diagonal, it suffices to consider 
$$N^{1/6}\le |u|\le \ep N^{1/2}\,.$$
As in the proof of Theorem \ref{neardiag}, we use a smooth partition of unity to
write  (\ref{SN}) as the sum of two integrals $I_1,I_2$.  Note that in the
proof that
$I_2=O(N^{-\infty})$, we use only that $|u|\le \ep N^{1/2}$. So it remains to
consider the integral $I_1$.  
Recalling (\ref{entirephase}), we have
$$\wt\Psi = \al_N(\theta)t-\theta\;,\qquad|\al_N(\theta)|
\ge \frac{|u|^2}{3N} \ge \frac{1}{3}N^{-2/3}\qquad \mbox{for}\ \ N^{1/6}\le |u|\le
\ep N^{1/2}\;.$$ Integrating by parts, we
then obtain
\begin{eqnarray*} &&\left|\int_0^3 
 e^{ i N \wt\Psi} \rho_1(t)s\big( \frac{u}{\sqrt{N}}, \theta; 
\frac{v}{\sqrt{N}}, 0, Nt\big) dt\right| \\ && \qquad\le
|N\al_N|^{-k} \sup_{t\le 3}
\left|\frac{\d^k}{\d t^k}\left[
\rho_1(t)s\big(
\frac{u}{\sqrt{N}},
\theta; 
\frac{v}{\sqrt{N}}, 0, Nt\big)\right]\right| \le C_k
N^{-k/3+m}\,,\end{eqnarray*} Hence by  (\ref{I1}),
\begin{equation}\label{outsideball} |I_1|\le C'_k N^{-k/3+m+1}\;.\end{equation}
and thus
$S_N^p(\frac{u}{\sqrtn},0;0,0) = O(N^{-\infty})$.  The estimate for the
derivatives is similar.\end{proof}

\begin{rem} It is of some interest to know how precise these off-diagonal
estimates  actually are. In the holomorphic case (with a $\ccal^\infty$ but not
necessarily analytic metric), M. Christ has shown (cf.\ \cite{Ch}) that one has
the global off-diagonal bounds
\begin{equation} \label{CHRIST} |\Pi_N(x,y)| \leq C N^m e^{- \sqrt{N} d(x,y)},
\;\;\;\;\; ((x,y) \in X \times X) \end{equation} for a constant $C$ which is
uniform in $(N, x, y).$ He also shows that this estimate cannot be improved in
general. Thus, the decay estimates of Lemma \ref{DECAY} do not hold in general
without the $O(N^{-\infty})$ term. It is quite possible that  (\ref{CHRIST})
extends to the almost complex context of this paper, since it is  based on
weighted $\lcal^2$ estimates for an associated Green's kernel which seem to have
natural analogues for the pseudodifferential
$\bar{D}_j$-complex of
\cite{BG}.  However, we do not need such a precise estimate as (\ref{CHRIST}) and
hence will not investigate it further.

\end{rem}

We now outline how Proposition \ref{transversal} follows from Lemma \ref{DECAY}.
The main steps in adapting the proof of Donaldson's
transversality theorem
\cite[Theorem~5]{DON.1} are to show the following:  Let
$$\textstyle\gcal_{p_i}^N(q) = |\sigma_{p_i}^N(q)+
|\dbar\sigma_{p_i}^N(q)|+\frac{1}{\sqrtn}|\nabla\sigma_{p_i}^N(q)|+
\frac{1}{\sqrtn} |\nabla\dbar\sigma_{p_i}^N(q)|
+\frac{1}{N}|\nabla^2\sigma_{p_i}^N(q)|
\;;$$ then
\begin{itemize}

\item  $\sum_i\gcal_{p_i}^N(q)$ is bounded independently of $N$;

\item  for all $D>0$ there exists $N_0(D)$ such that $\gcal_{p_i}^N(q) \leq
Ce^{-  D^2/3}$ for
$d_N(q,p_i)
\geq D$, $N>N_0(D)$.
\end{itemize}

These statements follow immediately from  parts (iii) and (iv) of Lemma
\ref{DECAY}; to obtain  the first statement we also note that
the number of points of the `lattice'
$\Lambda_N = \{p_i\}$ is $\sim N^m$ and that we can choose $D>0$ so that the balls
$\{B(p_i,\frac{D}{\sqrtn})\}$ cover $M$ and at most
$(2D+2)^{2m}$ of them pass through any point. 
{From} there (as in  the holomorphic case
described in \cite[p.~697]{DON.1}), we obtain Proposition \ref{transversal} from the
argument in \cite[\S\S3--4]{DON.1}.\qed

\subsection {The Adjunction Formula} \label{s-genus}

As explained in \cite{DON.1},
the zero set $\zb(s_N)$  of a section $s_N\in\ccal^\infty(M,L^N)$ satisfying the
transversality condition (\ref{QT}) is a symplectically embedded
submanifold.  In particular, when  
$(M, \omega)$ is a 4-manifold and
$\zb(s_N)$ is a symplectically embedded Riemann surface, then it is well known
that the genus $g_N$ of 
$\zb(s_N)$ is given by the Adjunction Formula as in the
holomorphic case:

\begin{equation} \label{genus}  g_N=\frac{c_1(L)^2}{2}N^2 -\frac{c_1(M)\cdot
c_1(L)}{2}N+1\,.
\end{equation}
To verify (\ref{genus}), we let $C$ denote the complex curve $\zb(s_N)$ with complex
structure
$J_0$ induced by
$\om$.   Then $J_0$ can be extended to
an almost complex structure $J'$ on $M$, compatible with $\om$, and it follows 
that
$C$ is a J-holomorphic curve in $(M,\om,J')$.  By our assumptions, the
homomorphism
$\nabla s_N:T_M\to L^N$ is surjective along $C$.  Hence, for all  $z_0\in C$,
 we have$$\ker_{z_0} \Big (\d s_N:T^{1,0}_M\to L_N\Big) =
\left(\ker_{z_0}\nabla s_N\right)  \cap T^{1,0}_{M,z_0} = T_{C,z_0} \cap
T^{1,0}_{M,z_0} =T^{1,0}_{C,z_0}$$  (where we use the almost complex structures
$J',J_0$ on
$M,C$, respectively). Proceeding as in the holomorphic case, from the resulting
exact sequence of vector bundles
$$0\to T^{1,0}_C \to  T^{1,0}_M|_C\, {\buildrel \d s_N \over\longrightarrow}\,
L^N|_C\to 0\,,$$ we obtain the Adjunction Formula:
\begin{equation} \label{adjunction} (K_M\otimes L^N)|_{C} =
K_{C}\,.
\end{equation}
Here $K_M=\bigwedge^2 \left.T^{1,0}_M\right.^*$ and $K_{C}=
\left.T^{1,0}_C\right.^*$.  Note that $c_1(K_M)$ is  independent of the choice of
$J'$, since the space of compatible almost complex structures on $(M,\om)$ is
connected (and in fact is contractible; see \cite[4.1]{McSa}). Thus we
have
$c_1(K_M)=-c_1(M)$. It follows further from Donaldson's Lefschetz hyperplane theorem  for
asymptotically holomorphic sections satisfying (\ref{QT}) (\cite[Prop.~39]{DON.1}) that $C$ is
connected.  We then conclude from (\ref{adjunction}) that 
$$2g-2=c_1(K_C)= \big(-c_1(M)+Nc_1(L)\big)\cdot Nc_1(L)\,,$$ 
which yields (\ref{genus}). 

\medskip
Auroux
\cite{A} generalized Donaldson's construction to obtain asymptotically holomorphic,
 `quantitatively transversal' sections of
$E\otimes L^N$, where E is an arbitrary complex vector bundle of rank $k\le m$ over
$M$. A section
$s_N\in\ccal^\infty(M,E\otimes L^N)$ is quantitatively transversal (\cite{A},
Definition~2) if $\nabla s_N:T_M\to E\otimes L^N$ is surjective and has a right
inverse of norm $<\frac{1}{\eta\sqrtn}$ at all points where $|s_N|<\eta$.
The zero sets $\zb(s_N)$ are then embedded symplectic submanifolds of $M$
(\cite{A}, Prop.~1). In particular, for the case
where $E$ is the trivial vector bundle $M\times \C^{k}$, sections of $E\otimes
L^N$ are $k$-tuples of sections of $L^N$, and we can
write $s_N=(s_N^1,\dots,s_N^k)$. Then $\zb(s_N)$ is the set of simultaneous
zeros of the sections $s_N^j\in H^0_J(M,L^N)$ and $s_N$ is asymptotically
holomorphic if the  $s_N^j$ are asymptotically holomorphic (but quantitative
transversality of the  $s^j_N$ is not sufficient to guarantee transversality of
$s_N$).  

By the previous comments, Auroux's argument can be used to give
asymptotically holomorphic, quantitatively transverse sections $s_N\in
H^0_J(M,L^N)^k$.  In the case where $k=m-1$, $\zb(s_N)$ is a complex curve, which
is connected by Auroux's generalization \cite[Prop.~2]{A} of  Donaldson's Lefschetz
hyperplane section theorem cited above, and the genus  formula (\ref{genus})
extends to this case as follows:

\begin{prop} {\rm (see \cite[Proposition 5]{A})} Let $s^1_N,\dots,s^{m-1}_N\in
\ccal^\infty(M,L^N)$ such that
$s_N:=(s^1_N,\dots,s^{m-1}_N)$ is quantitatively transverse and asymptotically
holomorphic, for $N=1,2,3,...$.  Then for
$N$ sufficiently large,
$\zb(s_N)$ is connected and its genus $g_N$ is given by
$$ g_N=\frac{m-1}{2} c_1(L)^m N^m -\half{c_1(M)\cdot
c_1(L)^{m-1}}N^{m-1}+1\,.$$ \end{prop}

The proposition follows by the argument in \cite{A}, which is a generalization of
the proof of (\ref{genus}).  To summarize this argument, the hypotheses imply that
$\zb(s_N)$ is a symplectic submanifold (for $N\gg 0$) and one obtains as before an
exact sequence
\begin{equation}\label{exact} 0\to T^{1,0}_C \to  T^{1,0}_M|_C\, {\buildrel \d s_N
\over\longrightarrow}\, L^N\oplus\cdots\oplus L^N|_C\to 0\,,\end{equation}
which yields
\begin{eqnarray*}2-2g_N & =& c_1(C)\ =\ \Big( c_1(M)-(m-1)Nc_1(L),[C]\Big)\\
&=&\big[ c_1(M)-(m-1)Nc_1(L)\big]\cdot c_1(L)^{m-1}N^{m-1}\,.\end{eqnarray*}

\begin{rem} The results of this
paper can be generalized to the case of almost holomorphic 
sections in
$H^0_J(M,E\otimes L^N)$, for an arbitrary complex vector bundle $E\to M$. Auroux's
argument can then be used to find quantitatively  transversal, asymptotically
holomorphic sections
$s_N\in H^0_J(M,E\otimes L^N)$. The zero sets $\zb(s_N)$ are then
$(m-k-1)$-connected symplectic submanifolds, for $N\gg 0$; if $k=\mbox{rank}\,
E=m-1$, they are connected Riemann surfaces.  In this case, the genus formula
becomes:
\begin{equation} \label{genus-E} 2g_N-2 = \big[(m-1)Nc_1(L)+c_1(E)-
c_1(M)\big]\sum_{j=0}^{m-1} c_{m-1-j}(E)c_1(L)^jN^j\,.\end{equation}
To verify (\ref{genus-E}), we proceed as in \cite{A}, noting first that we have an
exact sequence
\begin{equation}\label{exact-E} 0\to T^{1,0}_C \to  T^{1,0}_M|_C\, {\buildrel \d s_N
\over\longrightarrow}\, E\otimes L^N|_C\to 0\,.\end{equation}
The genus formula (\ref{genus-E}) then follows from (\ref{exact-E}) by a routine
Chern class calculation.
\end{rem}

\end{document}